\documentclass{article} 
\usepackage{amsmath}
\usepackage{amssymb}
\usepackage{amsthm}
\usepackage{paralist}
\usepackage{geometry}		            
\usepackage{graphics, graphicx}
\usepackage{epsfig} 
\usepackage{epstopdf} 
\usepackage[colorlinks=true]{hyperref}
\hypersetup{urlcolor=blue, citecolor=red}
\allowdisplaybreaks
\usepackage{cleveref}
\usepackage{subcaption}
\usepackage{authblk}                    
\usepackage[backend=biber,style=numeric,sorting=nty]{biblatex} 
\newtheorem{theorem}{Theorem}[section]
\newtheorem{corollary}{Corollary}

\newtheorem{lemma}[theorem]{Lemma}

\theoremstyle{definition}
\newtheorem{definition}[theorem]{Definition}

\newcommand{\norm}[1]{\left\lVert#1\right\rVert}
\newcommand{\abs}[1]{\lvert #1 \rvert}
\DeclareMathOperator{\Spec}{Spec}
\DeclareMathOperator{\diverg}{div}
\DeclareMathOperator{\TV}{TV}

\DeclareMathOperator*{\argmin}{arg\,min}
\DeclareMathOperator{\cm}{cm}
\DeclareMathOperator{\mm}{mm}
\DeclareMathOperator{\Pois}{Pois}

\title{On multiple scattering in Compton scattering tomography and its impact on fan-beam CT}

\author[1]{Lorenz Kuger}
\author[2]{Ga{\"e}l Rigaud}
\affil[1]{Friedrich-Alexander-Universität Erlangen-Nürnberg, Department Mathematics}
\affil[2]{Universität Stuttgart, Department Mathematics, Germany}
\date{May 2022}

\begin{document}
\maketitle
\begin{abstract}
 The recent development of energy-resolving scintillation crystals opens the way to new types of applications and imaging systems. In the context of computerized tomography (CT), it enables to use the energy as a dimension of information supplementing the source and detector positions. It is then crucial to relate the energy measurements to the properties of Compton scattering, the dominant interaction between photons and matter. An appropriate model of the spectral data leads to the concept of Compton scattering tomography (CST). Multiple-order scattering constitutes the major difficulty of CST. It is, in general, impossible to know how many times a photon was scattered before being measured. In the literature, this nature of the spectral data has often been eluded by considering only the first-order scattering in models of the spectral data. This consideration, however, does not represent the reality as second- and higher-order scattering are a substantial part of the spectral measurement. In this work, we propose to tackle this difficulty by an analysis of the spectral data in terms of modeling and mapping properties. Due to the complexity of the multiple order scattering, we model and study the second-order scattering and extend the results to the higher orders by conjecture. The study ends up with a general reconstruction strategy based on the variations of the spectral data which is illustrated by simulations on a joint CST-CT fan beam scanner. We further show how the method can be extended to high energetic polychromatic radiation sources.
\end{abstract}

\section{Introduction}\label{sec:Intro}
In conventional CT, an X-ray source emits ionizing radiation towards an object under study while detectors located outside the object measure variations in the intensity of the incoming photon flux. The classical mathematical model for CT \cite{Natterer01A, Natterer01B} involves the reconstruction of the attenuation coefficient $\mu_E$ through an inversion of the ray transform $\mathcal{P}$ in the following model
\begin{equation}\label{eqn:rayTransform}
	(\mathcal{P} \mu_E) (\theta, \mathbf{s}) := \int_{\mathbb{R}} \mu_E(\mathbf{s}+t\theta)\,\mathrm{d}t = \log I(\mathbf{s}) - \log(g_0(\mathbf{s},\mathbf{d},E)),
\end{equation}
where $I(\mathbf{s})$ is the initial intensity at the source $\mathbf{s}$ and $g_0$, called ballistic or primary radiation, is the radiation intensity measured at the detector $\mathbf{d}=s+T\theta$.
The quantity $\mu_E$ herein is a sum of contributions due to different physical processes occurring between photons and matter inside the medium under study. These include scattering events, photoelectric absorption and pair production (see section 2.3 in \cite{Leroy11}) and depend on the photon energy and the nature of the material. For energies typical in CT, the predominant effects are photoelectric absorption (denoted by $\mu^{\mathrm{PE}}(\mathbf{x};E))$) and Compton scattering (denoted by $\mu^{C}(\mathbf{x};E)$) cf. \cite{Stonestrom81}. Hence the attenuation coefficient is
\begin{equation}\label{eqn:approxLinAttOne}
	\mu_E(\mathbf{x}) = \mu^{\mathrm{PE}}(\mathbf{x}; E) + \mu^{\mathrm{C}}(\mathbf{x}; E) = \mu^{\mathrm{PE}}(\mathbf{x}; E) + \sigma^{\mathrm{C}}(E) n_e(\mathbf{x}),
\end{equation}
where $\sigma^{\mathrm{C}}(E)$ is the total cross section of Compton scattering at energy $E$ and $n_e(\mathbf{x})$ is the electron density of the material at the point $\mathbf{x}$. For energy levels larger than $100$ keV and in-vivo materials the Compton effect dominates in \cref{eqn:approxLinAttOne}, see \cite{Leroy11}. Hence, assuming large enough initial energies, we may neglect photoelectric absorption and take $\mu_E = \sigma^C(E) n_e$.

While scattering effects are a natural source of noise in CT, they can also be proved useful if the scattering process is correctly added into the mathematical imaging model. Indeed, when Compton scattered photons are observed at an energy resolving detector placed outside the specimen, they can be exploited in a reconstruction of the electron density map, the arising method is known as Compton scattering tomography (CST).
The relevant relation between measurable and imaged quantity is given by the kinematics and geometry of Compton scattering. When a photon scatters with initial energy $E_0$, the scattering angle $\omega$ and the energy after scattering $E_\omega^{E_0}$ obey the Compton formula
\begin{equation}\label{eqn:ComptonFormula}
	E_\omega^{E_0} = \frac{E_0}{1 + \frac{E_0}{m_e c^2}(1 - \cos \omega)},
\end{equation}
where $m_ec^2$ is an electron's energy at rest. This relates scattering geometry and photon energy, making it possible to deduce the potential trajectories and scattering points of a detected photon when its energy is measured. The fraction of photons that are scattered making an angle $\omega \in (0,\pi]$ with the incident direction is proportional to the electron density $n_e$ \cite{Klein29}.
Hence, we may use the energy measurements in a density reconstruction. In comparison to a classical CT architecture, a CST setup requires that the detectors be non-collimated and measure the energetic distribution of incoming photons. We speak of an energy spectrum measurement, denoted $\Spec$.

One of the major obstacles in CST is the presence of multiple scattering. 
The contributions of successive Compton scattering events have been studied in physics for Compton profile evaluation, but may not be applied directly to image formation in Compton scattering imaging, see for example \cite{PhysRevA.13.335,PhysRevA.14.313,PhysRevA.14.328}.
Many research works in the literature have modeled and analyzed the use of photons that are scattered once in the examined object \cite{Nguyen10, Norton94, Wang99, Truong07, Truong12, Webber15, Brateman84, Lale59, Clarke73, Palamodov11, Evans97, Rigaud17, Tarpau20a, Webber19,WebberIPI19, Webber18, Webber20}.  
This simplifies the reconstruction task, but is not a realistic setting as photons can be scattered more than once. When $g_0$ denotes ballistic radiation, the spectrum measured at a detector can be written as
\begin{equation}\label{eqn:DecompositionSpectrum}
	\Spec(\mathbf{s}, \mathbf{d}, E) = \sum_{i=0}^\infty g_i(\mathbf{s}, \mathbf{d}, E),
\end{equation}
where $g_i$ gives the flux density of all photons scattered exactly $i$ times and measured by the detector at $\mathbf{d}$. Multiple Compton scattering effects lead to more complex combinations of scattering sites and thus randomize the photons' energetic and spatial distribution. The implied intuition is thus that the different components $g_i$ in the spectrum deliver information about the object's electron density at different levels of smoothness. With increasing order of scattering, the total radiation intensity from $g_i$ spreads out and the smoothness increases, making the inversion more ill-posed, which will be made precise later. 

Following on from the recent work \cite{Rigaud21} which studied the problem of multiple scattering in 3D CST, the aim of this paper is to infer analytic properties of the contribution of multiply scattered photons to $\Spec(\mathbf{s},\mathbf{d},E)$ that can be of help in a reconstruction task in a two-dimensional setting. The advantage of our present study is fourfold: (i) the conditions over the smoothness properties for the first-order scattering are simpler than in 3D which simplifies the conception of future 2D CST-scanners; (ii) the general reconstruction strategy similar to what was proposed in \cite{Rigaud21}  is here illustrated on more realistic data since the energy resolution is assumed here constant and technologically feasible; (iii) while the source is first assumed monochromatic, i.e. all photons are emitted with the same energy $E_0$, the flexibility of the proposed modelings and reconstruction technique allows  an extension to polychromatic sources with sufficiently distinct characteristic peaks in the spectrum; and (iv) the combination of CT with CST helps to circumvent the non-linearity issue of the modelling by delivering a first approximation of the sought-for electron density.  

The manuscript is organized as follows. In \cref{sec:modelg1}, we recall the mathematical model of $g_1$ in two dimensions which can be represented as a weighted circular Radon transform, noted $\mathcal{T}_1$. To emphasize the problem of multiple scattering, we implement a standard minimization problem regularized by total-variation (TV) applied on $\Spec$ with the sole model of the first-order scattering. The results show severe artifacts and deformations in the reconstructed image. In order to tackle this difficulty, we derive the modeling, noted $\mathcal{T}_2$, of the second-order scattering (higher orders are much harder to model analytically and by nature purely stochastic) in  \Cref{subsec:forwardSecondOrder}, see \Cref{thm:secondOrderOperator}. An important obstacle in CST is the non-linearity of the forward operators. Therefore, we propose to study their smoothness properties in \cref{subsec:MappingProperties} via linearized versions noted here $\mathcal{L}_1^{\tilde{n}_e}$ and $\mathcal{L}_2^{\tilde{n}_e}$ respectively and defining Fourier integral operators (FIO). This approximation is made possible in practice by prior informations obtained here by standard CT. The mapping properties are demonstrated in \Cref{thm:sobolevtheoreml1,thm:sobolevtheoreml2} and given on an $L^2$-Sobolev scale. From the observation on the smoothness arises a heuristic reconstruction idea presented in \cref{subsec:Algorithm}. Following the intuition on the higher order terms above, we propose to apply a differential operator to forward model and data, cancelling out smoother terms and preserving valuable variations in $g_1$ that are not present in higher-order terms. By conjecture, the aforementioned smoothness properties would extend to the higher order terms $g_i, i\ge 3$, where we additionally note that due to the energy loss, it becomes less likely for photons to scatter again with increasing number of scattering events. Therefore, higher-order components $g_i, i\ge 3$ are naturally also taken care of in the final reconstruction algorithm. Furthermore, the inverse problem arising in CST is typically nonlinear due to attenuation effects. Instead of neglecting these effects and in order to use the measured data efficiently, we propose to also employ classical CT reconstruction (obtained through the component $g_0$). A combination of CT and CST in the same fan-beam scanning geometry with very small numbers of source positions and detectors is suggested and simulation results for various toy objects validate the derived approach for monochromatic and polychromatic sources in \cref{sec:Results}.

\section{The First-Order Scattering Model and Its Limitations}\label{sec:modelg1}

\begin{figure}[t]\centering
	\includegraphics[width=0.7\linewidth]{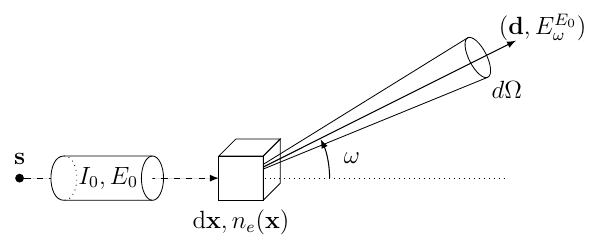}
	\caption{Illustration of the differential cross-section for the Compton effect.\label{fig:ComptonPhysics}}
\end{figure}

We are interested in modeling the flux of photons that are scattered once inside the specimen. The Compton total cross section $\sigma_c$ is a measure for the probability that a photon will undergo Compton scattering. In a thin layer of thickness $\mathrm{d}z$, the scattering probability is the fraction of intensity that is removed from the incident beam due to scattering:
\begin{equation}\label{eq:probaCompton}
	- \frac{\mathrm{d}I}{I} = n_e \sigma_c(E)\, \mathrm{d}z 
	\quad \Leftrightarrow \quad
	- \frac{\mathrm{d}I}{\mathrm{d}z} = n_e \sigma_c(E) I = \frac{\mathrm{d}n_s}{\mathrm{d}\mathbf{x}}
\end{equation}
where $\frac{\mathrm{d}n_s}{\mathrm{d}\mathbf{x}}$ is the number of photons scattered in the volume unit $\mathrm{d} \mathbf{x}$ and $n_e$ is the electron density. As depicted by Figure \ref{fig:ComptonPhysics}, the number of photons scattered from an incident beam of flux $I$ in the volume element $\mathrm{d}\mathbf{x}$ within a solid angle $\mathrm{d}\Omega$ oriented by $\omega$ is calculated by differentiating \cref{eq:probaCompton}.
\begin{equation}\label{eq:probaCompton_solidangle}
	\left(\frac{\mathrm{d}^2 n_s}{\mathrm{d}\mathbf{x}\, \mathrm{d}\Omega_c}\right)_\omega = I \cdot n_e \left( \frac{\mathrm{d} \sigma_c}{\mathrm{d}\Omega} \right)_\omega,
\end{equation}
with $\mathrm{d}\Omega = r^{-2} \mathrm{d}S$  a differential solid angle characterized by $r$, the distance of propagation and $\mathrm{d}S$, a differential element on the sphere. The last factor on the right hand side is the differential cross section $(\mathrm{d} \sigma_c / \mathrm{d}\Omega)_\omega = r_e^2 P(\omega)$ first described by Klein and Nishina \cite{Klein29}. $P(\omega)$ is sometimes also called Klein-Nishina probability.

We have to take into account the loss of intensity on the linear paths between source, scattering site and detector. When being emitted at $\mathbf{s}$ in direction $\mathbf{x}-\mathbf{s}$, the flux density is reduced due to photometric dispersion and the linear attenuation on the path between $\mathbf{s}$ and $\mathbf{x}$. 
This is modeled by multiplying the initial intensity by the factor
\begin{equation}\label{eqn:physicalFactors}
	A_{E_0}(\mathbf{s},\mathbf{x}) := \norm{\mathbf{s}-\mathbf{x}}^{-2}\exp\left( -\norm{\mathbf{s}-\mathbf{x}} \int_0^1\mu_{E_0}(\mathbf{s} + t(\mathbf{x}-\mathbf{s})) \,\mathrm{d} t \right).
\end{equation}

\begin{figure}[t]
	\centering
	\begin{subfigure}[b]{0.5\linewidth}
		\includegraphics[width=\linewidth]{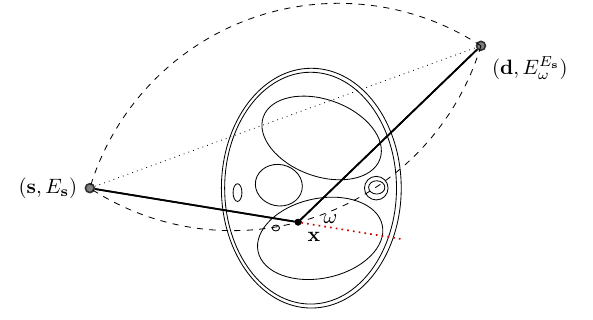}
		\caption{}
		\label{subfig:geometryFirstScattering}
	\end{subfigure}
	\hspace{0.3cm}
	\begin{subfigure}[b]{0.25\linewidth}
		\includegraphics[width = \linewidth]{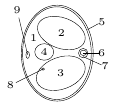}
		\caption{}
		\label{subfig:ThoraxIntensities}
	\end{subfigure}\\
	\vspace{0.2cm}
	\begin{tabular}{c|ccccccc} 
		\hline
		ellipse number & 1 & 2, 3 & 4 & 5 & 6 & 7, 9 & 8\\
		\hline
		$n_e/n_e^{\mathrm{water}}$ & 0.907 & 0.380 & 1.190 & 1.300 & 1.116 & 1.784 & 1.077 \\
		\hline
	\end{tabular}
	\caption{Geometry of the first order scattering (\cref{subfig:geometryFirstScattering}) for a human thorax phantom (\cref{subfig:ThoraxIntensities}).}
	\label{fig:firstOrderScGeometry}
\end{figure}

Combining eqs. (\ref{eq:probaCompton_solidangle}) and (\ref{eqn:physicalFactors}) leads to model the variation of the number of photons $g_1$ scattered at $\mathbf{x}$ and detected at $\mathbf{d}$ with energy $E_0$ by
\begin{equation}\label{eqn:probabilityPhotonScattersOnce}
	\frac{ \mathrm{d}^2 g_1(\mathbf{x}, \mathbf{d}, \mathbf{s})}{\mathrm{d}\mathbf{x} \,\mathrm{d}S} = I(\mathbf{s}) \underbrace{\left( \frac{\mathrm{d}\sigma}{\mathrm{d}\Omega} \right)_{\omega}  A_{E_0}(\mathbf{s} , \mathbf{x})A_{E_1}(\mathbf{x}, \mathbf{d})}_{=:w_1(n_e)(\mathbf{x},\mathbf{d},\mathbf{s})} n_e(\mathbf{x}).
\end{equation} 
The Compton formula \cref{eqn:ComptonFormula} allows a representation of the scattering angle $\omega$ in terms of the energies. For fixed positions $\mathbf{s}$ and $\mathbf{d}$ and a measured energy $E$, the set of possible scattering points are all $\mathbf{x}$ for which $\measuredangle(\mathbf{x}-\mathbf{s}, \mathbf{d}-\mathbf{x}) = \omega_{E}$ (where $\measuredangle(\mathbf{x}-\mathbf{s}, \mathbf{d}-\mathbf{x})$ is the angle between $\mathbf{x}-\mathbf{s}$ and $\mathbf{d} - \mathbf{x}$). Hence the locus of possible scattering events $\mathfrak{C}(E, \mathbf{d}, \mathbf{s})$ (dashed circular arcs in \cref{fig:firstOrderScGeometry}) is characterized by
\begin{equation*}
	\mathfrak{C}(E, \mathbf{d}, \mathbf{s}) = \left\{ \mathbf{x} \in \Omega\,:\, \frac{\mathbf{x}-\mathbf{s}}{\norm{\mathbf{x}-\mathbf{s}}}\cdot\frac{\mathbf{d}-\mathbf{x}}{\norm{\mathbf{d}-\mathbf{x}}} = \cos(\omega_E) = 1 - mc^2 ( 1/E - 1/E_0) \right\}.
\end{equation*}
The expected number of photons that are measured at $\mathbf{d}$ with energy $E$ can then be computed as the integral over $\mathfrak{C}$, weighted by the probability for a scattering event at each point.
By integrating \cref{eqn:probabilityPhotonScattersOnce} over all possible scattering sites $\mathbf{x} \in \mathfrak{C}$ that yield the same measured energy $E$, we obtain that
\begin{equation}\label{eqn:defT1}
	g_1(E,\mathbf{s},\mathbf{d}) \sim \mathcal{T}_1 (n_e) (E, \mathbf{s}, \mathbf{d}) := \int_{\mathfrak{C}(E, \mathbf{s},\mathbf{d})} w_1(n_e)(\mathbf{x},\mathbf{d},\mathbf{s}) n_e(\mathbf{x})\,\mathrm{d}\mathbf{x}
\end{equation}
where the weight function $w_1(n_e)(\mathbf{x},\mathbf{d},\mathbf{s}), x \in \mathfrak{C}(E,\mathbf{s},\mathbf{d})$ gathers the dispersion and attenuation factors in the model. Note that, by the definition in \cref{eqn:physicalFactors}, both $A_{E_0}$ and $A_{E_1}$ depend on the linear attenuation coefficient and hence on the electron density. The operator $\mathcal{T}_1$ is therefore nonlinear in $n_e$.

To obtain another, more convenient representation of the operator in order to be able later to analyze its mapping properties, we introduce a suited change of variable similarly as in \cite{Rigaud18}. The circular arcs are therein parametrized by $\mathbf{x} \in \mathfrak{C}(\omega, \mathbf{s}, \mathbf{d})\ \Leftrightarrow\ \cot \omega = \phi(\mathbf{x}-\mathbf{s},\mathbf{d}-\mathbf{s})$
where the characteristic function $\phi$ is given by
\begin{gather}\label{eqn:defphi}
	\phi(a,b) := \frac{\kappa(a,b) - \rho(a,b)}{\sqrt{1-\kappa(a,b)^2}},\quad\textrm{with }\kappa(a,b) := \frac{a^T b}{\norm{a}\norm{b}},\ \rho(a,b) := \frac{\norm{a}}{\norm{b}}.
\end{gather}
By the properties of the Dirac delta distribution $\delta$, we can reformulate \cref{eqn:defT1} to obtain
\begin{equation}\label{eqn:defT1delta}
	\mathcal{T}_1 (n_e) (\omega, \mathbf{s}, \mathbf{d}) = \int_{\Omega} \mathcal{W}_1(n_e)(\mathbf{x}, \mathbf{d}, \mathbf{s}) n_e(\mathbf{x}) \delta(\cot \omega - \phi(\mathbf{x}-\mathbf{s}, \mathbf{d}-\mathbf{s}))\,\mathrm{d} \mathbf{x},
\end{equation}
where the weight function depending on the electron density is now given by $$\mathcal{W}_1(n_e)(\mathbf{x}, \mathbf{d}, \mathbf{s}) := \norm{\nabla_{\mathbf{x}} \phi(\mathbf{x}, \mathbf{d}, \mathbf{s})} w_1(n_e)(\mathbf{x}, \mathbf{d}, \mathbf{s}).$$ 

Having established the integral operator $\mathcal{T}_1$, we wish to test how well it is suited as a forward model to reconstruct the electron density $n_e$ from the quantity $\Spec(\mathbf{s},\mathbf{d},E)$ that is expected to be measured in a realistic setting. As rough approximations to the energy spectrum, we test the model once with data $g \sim \Pois(g_1(\mathbf{s},\mathbf{d},E))$ (which, depending on the noise level, should give fairly good results with $\mathcal{T}_1$ as the forward model) and once with $g \sim \Pois((g_1 + g_2)(\mathbf{s},\mathbf{d},E))$ (which is closer to the realistic setting). Poisson distributed data $Y \sim \Pois(\lambda)$ is defined by the mass function
$$\mathbb{P}(Y = k) = \lambda^k e^{-\lambda}/k!$$
with respect to the counting measure on $\mathbb{N}_0$. The presence of Poisson noise requires some form of regularization. We add total variation (TV) regularization \cite{Rudin92, Burger13, Giusti84} and solve the variational problem
	\begin{equation}\label{eqn:variationalProblemOnlyT1}
		n_e^{\lambda} = \argmin_{n_e} \mathrm{R}(\mathcal{T}_1(n_e),g) + \lambda \TV(n_e),
\end{equation}
where $g$ is drawn from either $\Pois(g_1)$ or $\Pois(g_1+g_2)$. $\mathrm{R}$ is a suitable loss function. Since we assume a Poisson noise model here, we choose the Kullback-Leibler divergence 
	\begin{equation}\label{eq:kullbackLeiblerDivDef}
		\mathrm{R}(u,v) = \mathrm{KL}(u,v) = \int(u-v-v\log(u/v))\,\mathrm{d}x.
\end{equation}
The TV penalty is defined by
\begin{equation*}
	\TV(f) := \sup\left\{\int_\Omega f \diverg \phi\ \mathrm{d} x\ |\ \phi \in C_0^1(\Omega; \mathbb{R}^n):\abs{\phi(x)} \le 1 \textrm{ for } x \in \Omega \right\}.
\end{equation*}
We assume to measure a total of $8\cdot 10^6$ photons and compute synthetic data from \cref{eqn:defT1delta} for a phantom electron density $n_e$ modeling a transversal slice of a human thorax. The phantom (see \cref{subfig:ThoraxIntensities}) consists of characteristic functions of ellipses of different sizes and opacities and is a modified version of an earlier phantom which was used in \cite{Hahn14}. It is 28.4 cm resp. 21.3 cm wide at its largest and smallest diameters and its gray values are chosen as electron densities of materials typical in a human thorax \cite{Kanematsu16, xcomDatabase, Shrimpton81}.

\begin{figure}[t]
	\centering
	\begin{subfigure}{0.3\textwidth}
		\centering
		\includegraphics[width=\textwidth]{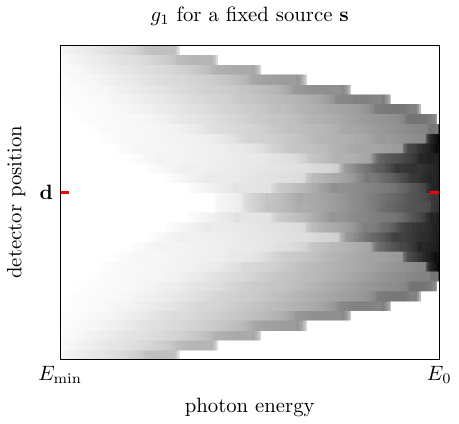}
		\subcaption{}
		\label{subfig:g1_singleSource}
	\end{subfigure}
	\hfill
	\begin{subfigure}{0.3\textwidth}
		\centering
		\includegraphics[width=\textwidth]{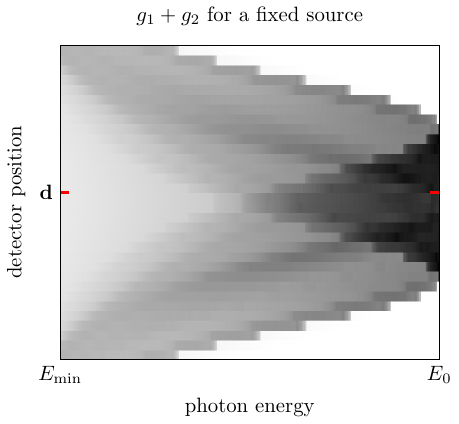}
		\subcaption{}
		\label{subfig:spec_singleSource}
	\end{subfigure}
	\hfill
	\begin{subfigure}{0.32\textwidth}
		\centering
		\includegraphics[width=\textwidth]{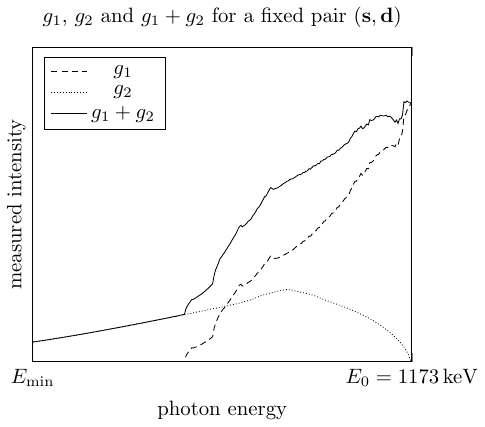}
		\subcaption{}
		\label{subfig:plots_data}
	\end{subfigure}\\
	\begin{subfigure}{0.31\textwidth}
		\includegraphics[width=\textwidth]{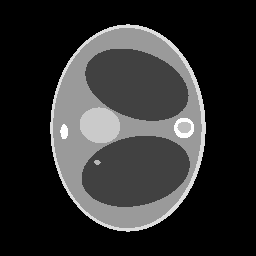}
		\subcaption{}
		\label{subfig:groundTruthThorax}
	\end{subfigure}
	\hfill
	\begin{subfigure}{0.31\textwidth}
		\includegraphics[width=\textwidth]{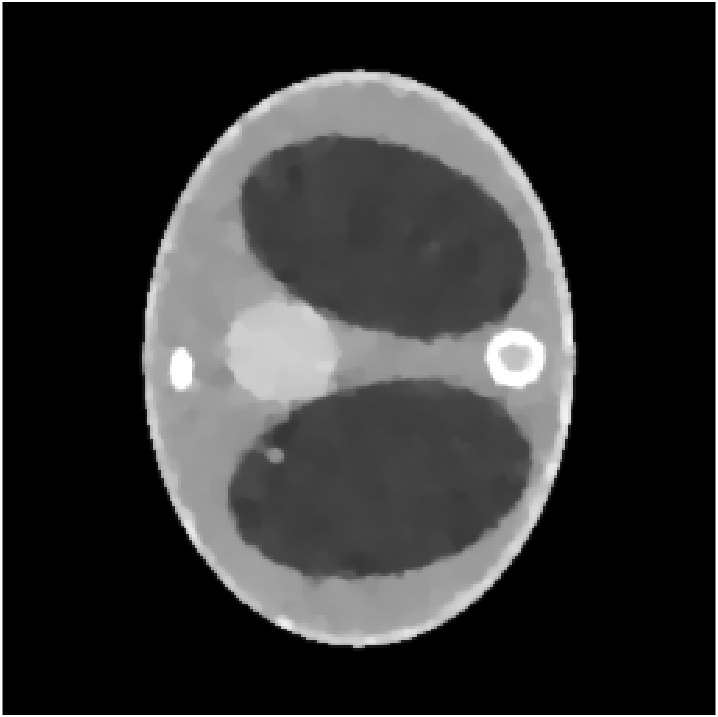}
		\subcaption{}
		\label{subfig:recOnlyT1_Poisg1}
	\end{subfigure}
	\hfill
	\begin{subfigure}{0.31\textwidth}
		\includegraphics[width=\textwidth]{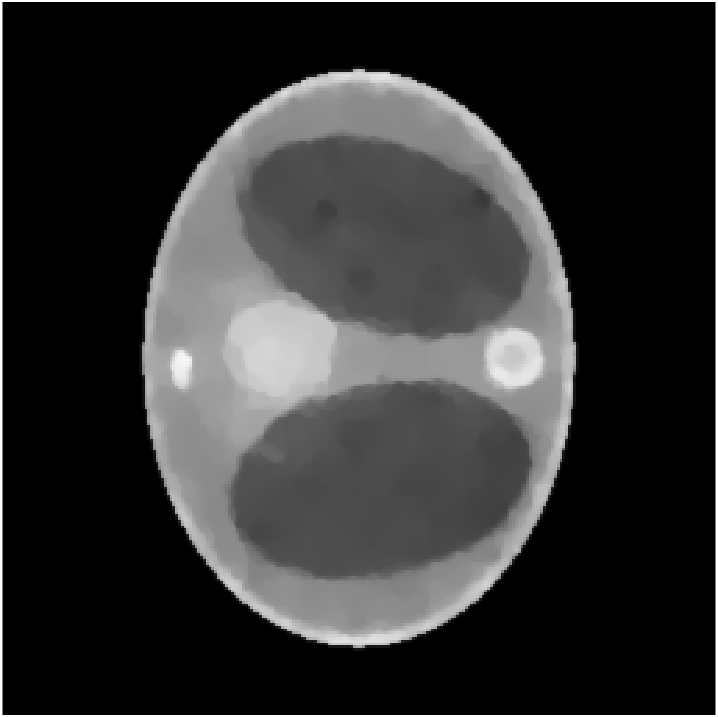}
		\subcaption{}
		\label{subfig:recOnlyT1_Poisg1g2}
	\end{subfigure}
	\caption{Top: Illustrations of the data $g_1$ in \ref{subfig:g1_singleSource} and $g_1+g_2$ in \ref{subfig:spec_singleSource}. $g_2$ alters the structure of the measured spectrum by adding large components not seen in $g_1$. The plot in \ref{subfig:plots_data} shows the measured spectrum for a fixed detector, corresponding to the row designated by $\mathbf{d}$ in the left two illustrations. Bottom: Reconstructions of the electron density of a thorax phantom (ground truth in \ref{subfig:groundTruthThorax}) using as forward model only $\mathcal{T}_1$ and data $\Pois(g_1)$ (\ref{subfig:recOnlyT1_Poisg1}) or $\Pois(g_1+g_2)$ (\ref{subfig:recOnlyT1_Poisg1g2}). Note how details of small size or lower contrast like the spine are badly imaged in the second setting.}
	\label{fig:reconstructionsUsingT1}
\end{figure}

In the reconstruction, a suitable choice for the regularization parameter $\lambda$ is computed by the L-curve method. In \cref{subfig:recOnlyT1_Poisg1} we see that a reconstruction of the electron density from $g = \Pois(g_1)$ works well and the TV penalty term reduces the effects from the Poisson noise sufficiently. However, adding the component $g_2$ (\cref{subfig:recOnlyT1_Poisg1g2}) to the spectrum distorts the reconstruction. Due to the enormous noise level brought into the model by the component $g_2$, the TV-regularized solution suffers from bad quality. Details are less visible or harder to localize and intensities of the different regions are altered complicating material recognition.

The second reconstruction shows how an analysis of the multiply scattered terms is inevitable if we wish to apply CST with data from a realistic setting.
In order to answer the question
\begin{center}\it 
	How to reconstruct the image in spite of double Compton scattering events?
\end{center}
we now derive an integral representation for $g_2$ and then analyze the smoothing properties of the components $g_1$ and $g_2$.

\section{Smoothness Properties of the Energy Spectral Data}\label{sec:smoothnessProperties}

The construction of the operator $\mathcal{T}_2$ modeling $g_2$ follows a comparable strategy as the one for $\mathcal{T}_1$. We identify possible pairs of scattering points and compute the integral weighted by the probability to detect a photon that is scattered at these two points. 

\subsection{Model of the Second Order Scattering}\label{subsec:forwardSecondOrder}
The main tool for the former step is a double application of the Compton formula \cref{eqn:ComptonFormula} which relates the two scattering angles and hence the two scattering points. Assume that a photon is emitted at the fixed point source $\mathbf{s}$ with initial energy $E_0$, scattered first at $\mathbf{x}$, then at $\mathbf{y}$ and is finally measured at $\mathbf{d}$ with energy $E$ (see \cref{fig:SecondOrderScGeometry}). Then the photon energy after the first scattering event is $E_1 := E_{\omega_1}^{E_0}$ where $\omega_1$ is the first scattering angle. After the second scattering with angle $\omega_2$, the remaining energy has to be $E = E_{\omega_2}^{E_1}$. By rearranging, we obtain
\begin{equation}\label{eqn:sumCosines}
	\cos \omega_1 + \cos \omega_2 = 2 - m_e c^2 \left( \frac{1}{E} - \frac{1}{E_0} \right) =: \lambda(E) \in (0, 2),
\end{equation}
hence $\omega_2$ is given by $\omega_2(\omega_1) = \arccos(\lambda(E) - \cos(\omega_1))$. This allows us to express $\mathbf{y}$ in terms of $\mathbf{x}$, $\omega_1$ and the energy $E$ at $\mathbf{d}$, which we note down as a Lemma:

\begin{figure}[t]
	\centering
	\includegraphics[width=0.55\linewidth]{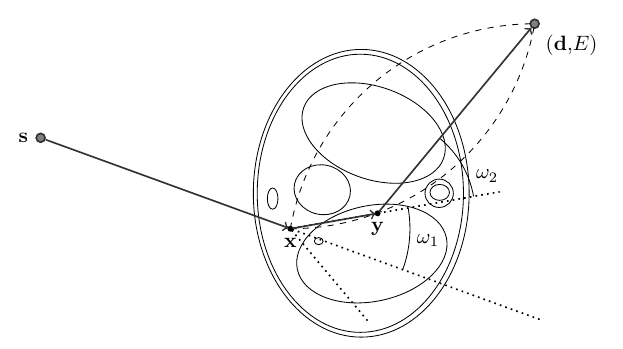}
	\caption{Geometry of the second order scattering}
	\label{fig:SecondOrderScGeometry}
\end{figure}

\begin{lemma}\label{lem:SecondScatteringPoint}
	Given the source $\mathbf{s}$, a detector $\mathbf{d}$, the measured photon energy $E$ and a first scattering point $\mathbf{x}$ with scattering angle $\omega_1$, we can compute the second scattering point by
	\begin{equation}\label{eqn:SecondScatteringPointY}
		\mathbf{y} = \mathbf{x} + r \begin{pmatrix} \sin ( \omega_1 - \beta + \pi/2)\\ \cos(\omega_1 - \beta + \pi/2)\end{pmatrix}\quad \textrm{if}\quad r > 0
	\end{equation}
	where the distance $r$ between $\mathbf{x}$ and $\mathbf{y}$ is given by
	\begin{equation}\label{eqn:SecondScatteringPointR}
		r = \norm{\mathbf{d} - \mathbf{x}} \left( \eta_2 - \cot( \omega_2(\omega_1)) \sqrt{1-\left(\eta_2\right)^2} \right).
	\end{equation}
	Herein, $\eta_2$ is the second component of $\eta := R^T (\sin ( \omega_1 - \beta + \frac{\pi}{2}),\cos(\omega_1-\beta+\frac{\pi}{2}))^T$, where
	$R$ is the orthogonal rotation matrix mapping $(0,1)^T$ to $(\mathbf{d} - \mathbf{x})/\norm{\mathbf{d} - \mathbf{x}}$ and $\beta$ is characterized by $\frac{\mathbf{x}-\mathbf{s}}{\norm{\mathbf{x}-\mathbf{s}}} = (\cos \beta, \sin \beta)^T$.
\end{lemma}

\begin{proof}
	With the first scattering parameters and the energy $E$ fixed, the second scattering point $\mathbf{y}$ has to lie on the 2D cone opening at $\mathbf{x}$ with angle $\omega_1$, but also on the circular arcs connecting $\mathbf{x}$ and $\mathbf{d}$ with opening angle $\omega_2(\omega_1)$. The calculations amount to computing the intersection point of these sets, details are deferred to \cref{pro:appxModelT2}.
\end{proof}

A special case mentioned only indirectly in \Cref{lem:SecondScatteringPoint} is when the distance between $\mathbf{x}$ and $\mathbf{y}$ is computed to be $r \le 0$, as the definition of $\mathbf{y}$ was only given for positive $r$. $r = 0$ actually corresponds to the case of single scattering and $r < 0$ represents the infeasible case where the photon loses too much energy already in the first scattering, i.e. $E_1 = E_{\omega_1}^{E_0} < E$, which can be discarded as no second scattering point can exist.
The proof implicitly uses $r > 0$ and for stability reasons, the numerical model will assume that no second scattering can occur in a small neighbourhood $\mathfrak{N}(\mathbf{x})$ of the first site, by just restricting $r > \varepsilon$ with a small $\varepsilon > 0$. We can now give a representation for the expected number of photons measured at the detector site using the previously computed relations:

\begin{theorem}\label{thm:secondOrderOperator}
	Assume we have a scanning architecture with a monochromatic, isotropic point source $\mathbf{s}$ emitting a fixed number of photons, a point detector at $\mathbf{d}$ and the electron density $n_e(\mathbf{x})$ is compactly supported on $\Omega \subset \mathbb{R}^2$. 
	Then the expected number of detected photons that are scattered exactly twice before being measured with an energy $E$, i.e. $g_2(E,\mathbf{d},\mathbf{s})$, is proportional to
	\begin{equation}\label{eqn:defS}
		\mathcal{T}_2 (n_e) (E,\mathbf{d},\mathbf{s}) = \int_{\Omega} \int_{\Lambda(E)} w_2(\mathbf{x}, \omega_1;\mathbf{d}, \mathbf{s}) n_e(\mathbf{x}) n_e(\mathbf{y}(\omega_1)) \,\mathrm{d} l(\omega_1)\,\mathrm{d}\mathbf{x}.
	\end{equation}
	The integration domain of $\omega_1$ is the set of possible first scattering angles $\Lambda(E) := \{\omega_1\in [-\pi, \pi]\,:\,E^{E_0}_{\omega_1} < E\}$, the weight function is 
	\begin{equation*}
		w_2(\mathbf{x}, \omega_1; \mathbf{d}, \mathbf{s}) = P(\omega_1)P(\omega_2(\omega_1))A_{E_0}(\mathbf{s},\mathbf{x}) A_{E_1}(\mathbf{x}, \mathbf{y}(\omega_1)) A_{E}(\mathbf{y}(\omega_1), \mathbf{d})
	\end{equation*}
	and the differential line segment is $\mathrm{d} l(\omega_1) = \sqrt{r^2 + \left(\frac{\partial r}{\partial\omega_1}\right)^2}\,\mathrm{d}\omega_1$, $r$ as in \Cref{lem:SecondScatteringPoint}.
\end{theorem}

\begin{proof}
	See  \cref{pro:appxModelT2}.
\end{proof}

Finally, we want to derive an alternative representation of $\mathcal{T}_2$ relating a measured energy $E$ to an integral over the whole domain $\Omega^2$ using a suited phase function. In the case of $\mathcal{T}_2$, we have to define the integration domain more carefully. 
As remarked after \Cref{lem:SecondScatteringPoint}, we discard the trivial case $\mathbf{y} = \mathbf{x}$ as it belongs actually to the first order scattering. 
Furthermore, we want to omit all other cases in which the points $\mathbf{s}, \mathbf{x}$ and $\mathbf{y}$ lie on the same line.
These correspond to the degenerate cases where the first scattering angle is exactly $0$ or $\pi$. Define therefore
\begin{equation}\label{eqn:Omega2}
	\Omega_2 := \left\{ (\mathbf{x},\mathbf{y}) \in \Omega^2\ |\ \mathbf{y}-\mathbf{x} \neq \lambda (\mathbf{x}-\mathbf{s}) \textrm{ for any } \lambda \in \mathbb{R}  \right\}.
\end{equation}
In order to obtain the phase function, we represent the scattering angles in terms of the points $\mathbf{s}, \mathbf{x}, \mathbf{y}$ and $\mathbf{d}$. Revisiting the definitions from \cref{eqn:defphi}, it holds $\cos \omega_1 = \kappa(\mathbf{y}-\mathbf{x}, \mathbf{x}-\mathbf{s}) = \frac{(\mathbf{y}-\mathbf{x})^T(\mathbf{x}-\mathbf{s})}{\norm{\mathbf{y}-\mathbf{x}}\norm{\mathbf{x}-\mathbf{s}}}$
and $\omega_2 = \cot^{-1}(\phi(\mathbf{y}-\mathbf{x},\mathbf{d}-\mathbf{x}))$
with $\phi = (\kappa-\rho)/\sqrt{1-\kappa^2}$. Using the relation \eqref{eqn:sumCosines} between the scattering angles and the measured energy, we obtain that a combination of points $\mathbf{s},\mathbf{x},\mathbf{y},\mathbf{d}$ yields a fixed energy $E$ if and only if
\begin{equation*}
	\lambda(E) = \psi(\mathbf{y}-\mathbf{x},\mathbf{x}-\mathbf{s},\mathbf{d}-\mathbf{x}) := \kappa(\mathbf{y}-\mathbf{x}, \mathbf{x}-\mathbf{s}) + \cos\left(\cot^{-1}(\phi(\mathbf{y}-\mathbf{x},\mathbf{d}-\mathbf{x}))\right).
\end{equation*}
Integrating over the feasible set $\Omega_2$ and using the new characteristic function $\psi$, we can therefore rewrite $\mathcal{T}_2$ as
\begin{equation}
	\mathcal{T}_2 (n_e) (E,\mathbf{d},\mathbf{s}) = \int_{\Omega_2} \mathcal{W}_2 (n_e) (\mathbf{x},\mathbf{y},\mathbf{d},\mathbf{s}) n_e(\mathbf{x}) n_e(\mathbf{y}) \delta(\lambda(E) - \psi(\mathbf{y}-\mathbf{x},\mathbf{x}-\mathbf{s},\mathbf{d}-\mathbf{x})) \,\mathrm{d}\mathbf{x}\mathrm{d}\mathbf{y}
\end{equation}
with the weight function $\mathcal{W}_2(n_e)(\mathbf{x},\mathbf{y},\mathbf{d},\mathbf{s}) = \norm{\nabla_{\mathbf{x},\mathbf{y}}\psi}w_2(\mathbf{x},\omega_1;\mathbf{d},\mathbf{s})$. Having established the model, we want to analyze its mathematical properties and in particular derive the $L^2$-Sobolev mapping theorems that will shed light on the nature of the different components $g_i$ in the energy spectrum. This will lead to the core idea of the reconstruction of the electron density.

\subsection{Mapping Properties}\label{subsec:MappingProperties}
\begin{figure}
	\centering
	\includegraphics[width=\linewidth]{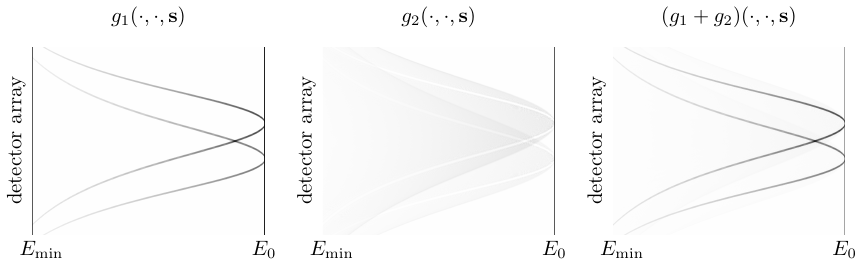}
	\caption{Contribution of $g_1$ and $g_2$ to the measured spectra of two small disks.}
	\label{fig:psfSpectra}
\end{figure}

As a motivational example for out method, consider an object $n_e^{\mathrm{test}}$ consisting only of two small disks in the plane.
The measured number of first and second order scattered photons can be simulated by evaluating the operators $\mathcal{T}_1$ and $\mathcal{T}_2$ of this object for a fixed source $\mathbf{s}$ and several detectors $\mathbf{d}$ and energies $E$. The data is visualized in \cref{fig:psfSpectra}.
Significant differences between $g_1 = \mathcal{T}_1 (n^{\mathrm{test}}_e)$ and $g_2 = \mathcal{T}_2 (n^{\mathrm{test}}_e)$ immediately become visible.
The spectrum of the second order scattered data is more '{}spread out'{} over the different energy levels and smoother in general.

The theory of Fourier integral operators (FIO) can be applied to study mapping properties of generalized Radon transforms. Expressing the kernel of the FIO in the form of an oscillatory integral allows to deduce $L_2$-Sobolev continuity results by simply examining the phase function, see e.g. \cite{Greenleaf2002}. In particular, the theory proved to be particularly useful in several imaging applications \cite{Krishnan15, Webber19, Hahn19, Kuchment95}. In this section, we apply it to the operators defining the different parts of the measured energy spectrum regarding multiple Compton scattering. The nonlinear operators $\mathcal{T}_1$ and $\mathcal{T}_2$ are first approximated by linear counterparts $\mathcal{L}_1, \mathcal{L}_2$. Instead of just omitting the nonlinearities due to attenuation effects entirely, it is more reasonable to approximate the nonlinear weight function through a prior estimate. A strategy on how the prior information can be obtained is explained in \cref{subsec:Algorithm}. The arising (linear) operators $\mathcal{L}_1$ and $\mathcal{L}_2$ are weighted Radon transforms and consequently constitute examples of FIOs. This allows to use the associated theory from \cite{Hoermander71, Krishnan15} and prove the $L_2$-Sobolev mapping properties. Our motivational idea from the beginning, which showed how $\mathcal{L}_2$ produces smoother data than $\mathcal{L}_1$ while the latter preserves the object's singularities better, reflects in these mapping properties.

\begin{definition}[Linearization of $\mathcal{T}_1$]
	Assume there is a smooth prior information $\tilde{n}_e\in C^{\infty}(\Omega)$ given which approximates the sought-for image.
	As a linear approximation of the first operator $\mathcal{T}_1$, we define
	\begin{equation}\label{eqn:defL1}
		\mathcal{L}_1^{\tilde{n}_e} (f) (p, \mathbf{s}, \mathbf{d}) = \int_{\Omega} \mathcal{W}_1(\tilde{n}_e)(\mathbf{x}, \mathbf{d}, \mathbf{s}) f(\mathbf{x}) \delta(p - \phi(\mathbf{x}-\mathbf{s}, \mathbf{d}-\mathbf{s})) \,\mathrm{d} \mathbf{x}
	\end{equation}
	where the only difference to $\mathcal{T}_1$ lies in the weight function $w_1 = \mathcal{W}_1(\tilde{n}_e)$ depending not on the function $f$, but rather on the prior.\\
	Note that the first argument is now $p = \cot \omega$ instead of the energy $E$. This simplifies notation and is not problematic since the three parameters $p, \omega$ and $E$ are all related by bijective mappings.
\end{definition}

\begin{definition}[Linearization of $\mathcal{T}_2$]
	Similarly, we define the linear operator $\mathcal{L}_2^{\tilde{n}_e}$ for a prior $\tilde{n}_e\in C^{\infty}(\Omega)$. To shorten the notation, gather $\mathbf{z} := (\mathbf{x},\mathbf{y}) \in \Omega_2 \subset \mathbb{R}^4$ and denote $\tilde{f}(\mathbf{z}) := f(\mathbf{x})f(\mathbf{y})$. This yields an operator linear in $\tilde{f}$
	\begin{equation*}
		\mathcal{L}_2^{\tilde{n}_e} (\tilde{f}) (\lambda, \mathbf{s}, \mathbf{d}) = \int_{\Omega_2} \mathcal{W}_2(\tilde{n}_e)(\mathbf{x}, \mathbf{y}, \mathbf{d}, \mathbf{s}) \tilde{f}(\mathbf{z}) \delta(\lambda - \psi(\mathbf{y}-\mathbf{x},\mathbf{x}-\mathbf{s},\mathbf{d}-\mathbf{x})) \,\mathrm{d}\mathbf{z}
	\end{equation*}
	where $\lambda = \lambda(E)$ as before.
\end{definition}

We first recall the relevant definitions for FIOs from \cite{Krishnan15, Hoermander71}. For a function $f:\mathbb{R}^n\to \mathbb{R}$, denote $\partial_x f(x) = \sum_{i=1}^n \frac{\partial f}{\partial x_i} \,\mathrm{d} x_i$. Further, we write $H^\nu_c(\Omega)$ for the space of all $u\in L^2(\Omega)$ compactly supported in $\Omega$ that satisfy
\begin{equation*}
	\norm{u}_{H^\nu} = \left(\int_\Omega (1+\norm{\xi}^2)^\nu \norm{\hat{u}(\xi)}^2 \,\mathrm{d}\xi \right)^{1/2} < \infty
\end{equation*}
where $\hat{u}$ is the Fourier transform of $u$ and $H^\nu_{\mathrm{loc}}(\Omega)$ for the space of all $u$ supported in $\Omega$ such that $\chi u \in H^\nu_c(\Omega)$ for all compactly supported $\chi \in C^\infty_c(\Omega)$.

\begin{definition}[Phase function, \cite{Krishnan15}]
	Let $Y \subset \mathbb{R}^m$ and $X \subset \mathbb{R}^n$ be open. Then the function $\Phi \in C^\infty(Y\times X \times \mathbb{R}\backslash\{0\})$ is called phase function if 
	\begin{itemize}
		\item $\Phi(y,x,\sigma)$ is positively homogeneous of order 1 in $\sigma$ and
		\item $(\partial_y \Phi, \partial_\sigma \Phi)$ and $(\partial_x \Phi, \partial_\sigma \Phi)$ never vanish for any $(y,x,\sigma) \in Y\times X \times \mathbb{R}\backslash\{0\}$.
	\end{itemize}
	The phase function is further called non-degenerate if $\partial_{y,x}(\partial_\sigma \Phi)$ do not vanish for $(y,x,\sigma) \in \Sigma_\Phi$ where
	\begin{equation*}
		\Sigma_\Phi := \left\{ (y,x,\sigma) \in Y\times X\times \mathbb{R}\backslash\{0\}\ |\ \partial_\sigma \Phi = 0 \right\}.
	\end{equation*}
\end{definition}

\begin{definition}[Fourier Integral Operator, \cite{Krishnan15}]
	The operator $\mathcal{T}$ is called a Fourier integral operator of order $s-\frac{n+m-2}{4}$ if 
	\begin{equation}\label{eqn:defGeneralFIO}
		\mathcal{T}f(y) = \int e^{i\Phi(y,x,\sigma)} p(y,x,\sigma) f(x) \,\mathrm{d} x \,\mathrm{d}\sigma,
	\end{equation}
	where $\Phi$ is a non-degenerate phase function and the following property is satisfied:\\
	For the function $p \in C^\infty(Y,X,\mathbb{R})$, called symbol of order $s$, it holds that, for every $K \subset Y \times X$ and all $\alpha \in \mathbb{N}_0, \beta \in \mathbb{N}_0^n, \gamma \in \mathbb{N}_0^m$, there exists a constant $C = C(K,\alpha,\beta,\gamma)$ such that
	\begin{equation*}
		\left|D_\sigma^\alpha D_x^\beta D_y^\gamma p(y, x, \sigma)\right| \le C(1+\|\sigma\|)^{s-\alpha} \text { for all } (y, x) \in K \text { and } \sigma \in \mathbb{R}.
	\end{equation*}
\end{definition}

Phase functions of a special, semiaffine form (applying to our cases) allow an analysis of mapping properties of the corresponding FIO.
\begin{lemma}[\cite{Hoermander71}]\label{lem:sobolevMappingFIO}
	Consider the special case where $y = (p,\theta) \in \mathbb{R} \times \Theta =: Y$ with $\Theta \subset \mathbb{R}^{m-1}$ open and the phase function is of the form 
	\begin{equation*}
		\Phi(p,\theta,x,\sigma) = - \sigma(p-\phi(x,\theta)).
	\end{equation*}
	If the vectors $\nabla_x \phi, \partial_{\theta_1}\nabla_x\phi, \dots, \partial_{\theta_{m-1}}\nabla_x\phi$ are linearly independent for all $(p,\theta,x) \in \mathbb{R}\times\Theta\times X$, then 
	\begin{equation*}
		\mathcal{T}: H_{c}^{\nu}(X) \rightarrow H_{\mathrm{loc}}^{\nu-k}(Y)
	\end{equation*}
	is a continuous operator.
\end{lemma}

Now we fit the weighted Radon transforms $\mathcal{L}_1$ and $\mathcal{L}_2$ into the above framework. Begin with $\mathcal{L}_1$ and rewrite it to fit the definition of an FIO. Using the Fourier transform of the Dirac delta, we obtain from eq. \eqref{eqn:defL1} the representation
\begin{equation*}
	\mathcal{L}_1^{\tilde{n}_e} (f) (p, \mathbf{s}, \mathbf{d}) = \frac{1}{2\pi} \int_{\Omega} \int_{\mathbb{R}} \mathcal{W}_1(\tilde{n}_e)(\mathbf{x}, \mathbf{d}, \mathbf{s}) n_e(\mathbf{x}) e^{-i\sigma(p - \phi(\mathbf{x}-\mathbf{s}, \mathbf{d}-\mathbf{s}))} \,\mathrm{d} \sigma \,\mathrm{d} \mathbf{x}.
\end{equation*}
In particular, the phase function $\Phi(p,\theta,\mathbf{x},\sigma) = - \sigma(p - \phi(\mathbf{x}-\mathbf{s}, \mathbf{d}(\theta)-\mathbf{s}))$ is of the special form of \Cref{lem:sobolevMappingFIO}, which shows that it is applicable to prove the Sobolev mapping property of $\mathcal{L}_1$ if all the necessary conditions on the phase function are satisfied.

\begin{theorem}\label{thm:fiol1}
	Assume that $\tilde{n}_e \in C^\infty(\Omega)$ and that the source $\mathbf{s} \notin \Omega$ is fixed . Then the operator $\mathcal{L}_1^{\tilde{n}_e}$ is an FIO of order $-\frac{1}{2}$.
\end{theorem}

\begin{proof}
	The proof consists in the verification of the different properties of the symbol and the phase function and is straightforward, the details are deferred to  \cref{pro:ThmL1}.
\end{proof}

Given that $\mathcal{L}_1$ is an FIO, it remains to check the necessary properties of the phase for the Sobolev mapping property. For this, we first prove a Lemma which relates the main condition of \Cref{lem:sobolevMappingFIO}, namely the linear independence of $\nabla_x \phi, \partial_{\theta}\nabla_x \phi$, to an equivalent condition on the measure geometry.

\begin{lemma}\label{lem:immersionl1}
	Let the detectors be defined by
	\begin{equation}\label{eqn:detectorPositions}
		\mathbf{d}(\theta) = \mathbf{s} + t(\theta) \begin{pmatrix}\cos\theta\\\sin\theta\end{pmatrix}
	\end{equation}
	where $\theta \in \Theta$, $\Theta$ some open interval. Then, if we write $\mathbf{x}-\mathbf{s} = (r\cos \xi, r\sin \xi)^T$ and the condition
	\begin{equation}\label{eqn:conditionImmersionl1}
		r \neq t(\theta) \cos (\theta - \xi) - t'(\theta) \sin(\theta - \xi)
	\end{equation}
	is satisfied for all $\mathbf{x} \in \Omega$ and every $\theta \in \Theta$, we have that $\det(\nabla_{\mathbf{x}}\phi, \partial_\theta \nabla_{\mathbf{x}}\phi) \neq 0$ for all $(\mathbf{x},\mathbf{d}) \in \Omega \times \mathbb{D}$.
\end{lemma}

\begin{proof}
	See \cref{pro:ThmL1}.
\end{proof}

Having established condition \eqref{eqn:conditionImmersionl1} for the Sobolev mapping property, one can directly check whether it is satisfied for a given scanner setup with specified source and detector positions. We carry out the computations for a standard fan-beam geometry known from CT that will be used in the numerical experiments.

\begin{theorem}\label{thm:sobolevtheoreml1}
	Let $\tilde{n}_e \in C^\infty(\Omega)$. Let the detectors $\mathbf{d}(\theta)$ be positioned on the circle $S^1$, likewise the source $\mathbf{s}$, e.g. at $(-1,0)^T$, and let the object be supported in the open disk of radius 1 around 0, i.e. $\Omega \subset B_1(0)$. Then the operator $$\mathcal{L}_1^{\tilde{n}_e} : H_c^\nu(\Omega) \to H_{\mathrm{loc}}^{\nu+1/2}(\mathbb{R} \times S^1)$$ is continuous.
\end{theorem}

\begin{proof}
	The positions of the detectors relative to the source are given by $\mathbf{d}(\theta) = \mathbf{s} + 2 \cos \theta \begin{pmatrix} \cos \theta\\ \sin \theta \end{pmatrix}$
	where $t(\theta) = 2\cos \theta$ and $\theta \in (-\pi/2,\pi/2)$. The preceding \Cref{lem:immersionl1,lem:sobolevMappingFIO} prove the theorem if eq. \eqref{eqn:conditionImmersionl1} is satisfied for all $\mathbf{x}-\mathbf{s} = (r\cos\xi, r\sin\xi)\in \Omega$. The condition reads
	\begin{equation*}
		r \neq t(\theta) \cos (\theta - \xi) - t'(\theta) \sin(\theta - \xi) = 2\cos\xi,
	\end{equation*}
	leading to $\mathbf{x} \notin S^1$, which is true due to the assumption $\Omega \subset B_1(0)$.
\end{proof}

This concludes our study of $\mathcal{L}_1$. We now want to derive similar results for the second operator $\mathcal{L}_2$. Start again by rewriting the integral transform in the form of an FIO:
\begin{equation*}
	\mathcal{L}_2^{\tilde{n}_e} (\tilde{f}) (\lambda, \mathbf{s}, \mathbf{d}) = \frac{1}{2\pi} \int_{\Omega_2} \mathcal{W}_2(\tilde{n}_e)(\mathbf{x}, \mathbf{y}, \mathbf{d}, \mathbf{s}) \tilde{f}(\mathbf{z}) e^{-i\sigma(\lambda - \psi(\mathbf{y}-\mathbf{x},\mathbf{x}-\mathbf{s},\mathbf{d}-\mathbf{x}))} \,\mathrm{d}\sigma\,\mathrm{d}\mathbf{z}
\end{equation*}
where we used again the Fourier representation of $\delta$. In the context of \Cref{lem:sobolevMappingFIO}, the phase function is given by $\Psi(\lambda,\theta,\mathbf{x},\mathbf{y},\sigma) = -\sigma(\lambda - \psi(\mathbf{y}-\mathbf{x},\mathbf{x}-\mathbf{s},\mathbf{d}(\theta)-\mathbf{x}))$. The next theorem establishes that the operator is indeed an FIO.

\begin{theorem}\label{thm:fiol2}
	Let $\tilde{n}_e \in C^\infty(\Omega)$, then $\mathcal{L}_2^{\tilde{n}_e}$ is an FIO of order $-1$.
\end{theorem}

\begin{proof}
	As for $\mathcal{L}_1^{\tilde{n}_e}$, the proof consists in verifying the properties of the symbol and the phase. This is straightforward, but rather technical, and therefore moved to  \cref{pro:ThmFioSob}.
\end{proof}

The corresponding Sobolev mapping property can now be proved for a given scanning geometry, we verify that it is valid for the same fan-beam setup that was defined in \Cref{thm:sobolevtheoreml1}. Note that $\mathcal{L}_2^{\tilde{n}_e}$ maps half a step more on the Sobolev space scale than $\mathcal{L}_1^{\tilde{n}_e}$ when applied to the same function. This leads to an important insight about the measured data, namely that $g_2$ is smoother than $g_1$.

\begin{theorem}\label{thm:sobolevtheoreml2}
	Let $\tilde{n}_e \in C^\infty(\Omega)$ and the measure geometry be defined as in \Cref{thm:sobolevtheoreml1} - the detectors $\mathbf{d}(\theta)$ are distributed on $S^1$, $\mathbf{s}$ is fixed at $(-1,0)$ and $\Omega \subset B_1(0)$. Then the operator $$\mathcal{L}_2^{\tilde{n}_e} : H_c^\nu(\Omega_2) \to H_{\mathrm{loc}}^{\nu+1}((0,2) \times S^1)$$  is continuous.
\end{theorem}

\begin{proof}
	See  \cref{pro:ThmFioSob}.
\end{proof}

This concludes our analytic study of the operators.
The important aspect for our objective of carrying out electron density reconstruction is the smoothness of the different components $g_1$, $g_2$ as the approximations of $\mathcal{T}_1$, $\mathcal{T}_2$ map to $H^{1/2}_{\mathrm{loc}}$ and $H^{1}_{\mathrm{loc}}$, respectively. Although we only examined $\mathcal{L}_1$ and $\mathcal{L}_2$ here, the general idea can be extended further, as the emitted photons can of course scatter more often than two times in the object. Modeling the process of $n$ scattering events would most probably be very complicated and proving suitable mapping properties for the arising operators even harder. Nevertheless, going out from our study of the first and second order here, one can expect that, in terms of the Sobolev mapping properties elaborated in \Cref{thm:sobolevtheoreml1,thm:sobolevtheoreml2}, the higher order parts of the spectrum are contained at least in $H^1(Y)$ where $Y\subset \mathbb{R} \times \Theta$.

We are going to make use of the smoothness properties in the next section, where we derive a method to recover the image from the spectrum.

\subsection{Image reconstruction exploiting smoothness properties}\label{subsec:Algorithm}
We now want to examine how we can exploit the results on $g_1, g_2$ in image reconstruction. Efforts have been made to reconstruct images from simulated approximations to $g_1$ by inversion-type formulae or backprojection of $\mathcal{T}_1$ \cite{Rigaud18, Tarpau20a, Cebeiro21, Rigaud17}. Attenuation effects are often neglected in order to preserve linearity of the operator or need to be refined parallel to the density reconstruction \cite{Goedecke2022}. Moreover, we showed in \cref{sec:modelg1}, models including only first order scattering and ignoring higher order terms \textbf{fail on realistic data that includes $g_i, i \ge 2$}. We wish to address both problems, i.e. nonlinearity and multiple scattering, by combining CT and CST in the same scanning geometry and by exploiting the analytical results on the spectrum from \cref{subsec:MappingProperties}.

\subsubsection{Linearization using a sparse-view CT prior}
As we saw, it is possible to attain a linear operator $\mathcal{L}_1^{\tilde{n}_e}$ when some prior information $\tilde{n}_e$ is available. In order to obtain this prior information, we make use of classical CT and the ballistic radiation $g_0$ in the measured spectrum. We will work with a scanning geometry with very few numbers of source and detector positions. This is possible because the same scanning geometry is used for CT and CST, where the latter benefits from the energy as additional dimension of information. The small number of source positions further has practical advantages like the resulting smaller scanning times and radiation exposure.

We try to solve the CT problem \eqref{eqn:rayTransform}, i.e. find $\mu$ from measurements $g_0$ where $\mathcal{P}\mu_{E_0} = -\log(g_0/I(\mathbf{s}))$. In our setting with few source and detector positions, the size of $g_0$ is much smaller than the size of $\mu_{E_0}$ the problem is called underdetermined. Solving such a system by an FBP method leads to artifacts and bad image quality. Instead, we regularize the solution by employing total variation regularization \cite{Rudin92}, as has been done before for sparse-view CT \cite{Sidky08, Zhu13}. This consists in solving the variational problem

	\begin{equation}\label{eqn:CTLimitedDataMinimization}
		\sigma^{\textrm{C}}(E_0)\tilde{n}_e^\lambda = \tilde{\mu}_{E_0}^\lambda = \argmin_{f} \mathrm{R}(\mathcal{P}f,\log(g_0/I(\mathbf{s}))) + \lambda \textrm{TV}(f).
	\end{equation}

We obtain a prior estimate $\tilde{\mu}_{E_0}$ of the attenuation map, which is typically of suboptimal quality due to the sparse data. Via the first equality, this leads to an estimate $\tilde{n}_e$ of the electron density. $\tilde{n}_e$ is plugged into the exponential terms of $\mathcal{T}_1$ so that the weight function $\mathcal{W}_1$ is fixed, leading to the linearized operator $\mathcal{L}_1^{\tilde{n}_e}$ that was analyzed in \cref{subsec:MappingProperties}. $\mathcal{L}_1^{\tilde{n}_e}$ might constitute a sufficient (inexact) forward model in the task of recovering $n_e$ from $g_1$ if $\tilde{n}_e$ is sufficiently close to the ground truth $n_e^{\textrm{ex}}$. Using variational TV regularization, the corresponding inverse problem can then be solved by the following optimization problem
\begin{equation}\label{eqn:CSTvariationalForm1}
	n_e^{\lambda} = \argmin_{f} \mathrm{R}(\mathcal{L}_1^{\tilde{n}_e}(f), g) + \lambda \TV(f).
\end{equation}
If the data in \eqref{eqn:CSTvariationalForm1} were $g = g_1 := \Pois(\mathcal{T}_1(n_e^{\mathrm{ex}}))$, we could hope for accurate solutions, but following the discussions in section \ref{sec:modelg1}, considering realistic cases with $g = g_1+g_2+\dots$ leads to a severe model inexactness.

\subsubsection{Reducing the impact of multiply scattered terms}
We want to reconstruct $n_e$ when the data $g$ contains the second scattering component $g_2 = \mathcal{T}_2n_e^{\mathrm{ex}}$. In section \ref{sec:modelg1} we saw that a simple reconstruction using only the forward model $\mathcal{T}_1$ is not sensible. An inversion formula for the operator $\mathcal{T}_2$, let alone $\mathcal{T}_1+\mathcal{T}_2$, is unknown, reasons for this are the complex structure of the phase function and the nonlinearity of $\mathcal{T}_2$. On the other hand, including $\mathcal{T}_2$ or even $\mathcal{L}_2^{\tilde{n}_e}$ into the forward model seems computationally cumbersome due to the additional integration variable.

We propose a treatment of the multiply scattered terms by preprocessing the measured spectrum in the energy variable through another operator $\mathcal{D}_E$. This transforms \cref{eqn:CSTvariationalForm1} into
\begin{equation}\label{eqn:CSTvariationalForm2}
	n_e^{\lambda} = \argmin_{f} \mathrm{R}(\mathcal{D}_E \mathcal{L}_1^{\tilde{n}_e}(f),\mathcal{D}_E g) + \lambda \TV(f),
\end{equation}
where now $g \sim \Pois((\mathcal{T}_1 + \mathcal{T}_2)(n_e^{\mathrm{ex}}))$. Since the data points $\mathcal{D}_E g$ do not follow Poisson distributions anymore, the data fidelity is set to an $L_2$ norm in the numerical tests. It could be interesting to further investigate the statistical nature of $\mathcal{D}_Eg$ and the effect of choosing other data fidelity terms, e.g. weighted $L_2$ norms \cite{Wang2006}.

$\mathcal{D}_E$ should not increase the ill-posedness or ill-conditioning of the problem too much and preserve important information about the object from $g_1$. At the same time, it should reduce as much as possible the impact of $g_2$ on the reconstruction.
Following the intuition presented in the remarks on \cref{fig:psfSpectra} and the analytical results in \cref{subsec:MappingProperties}, the higher-order term $g_2$ is smoother than $g_1$ in terms of a Sobolev space scale. Hence, a heuristic choice for $\mathcal{D}_E$ is a differential operator $\mathcal{D}_E$ in the energy variable.
In practice, the differentiation step itself is ill-posed and needs to be regularized. For a low pass filter function $F_\gamma$, $\gamma > 0$, define the operator $\mathcal{D}_E^\gamma$ by
\begin{equation*}
	\mathcal{D}_E^\gamma g^{\mathbf{s},\mathbf{d}} := \mathcal{F}^{-1} (i \zeta \cdot F_\gamma(\zeta) \cdot \mathcal{F}(g^{\mathbf{s},\mathbf{d}})(\zeta))
\end{equation*}
where $g^{\mathbf{s},\mathbf{d}}$ is the restriction of the data to a given source-detector pair and $\mathcal{F}$ defined by $\mathcal{F}(g^{\mathbf{s},\mathbf{d}})(\zeta) := (2\pi)^{-1/2} \int g^{\mathbf{s},\mathbf{d}}(E) \exp(-i\zeta E) \,\mathrm{d} E $ is the Fourier transform w.r.t. the energy variable. 

\begin{figure}[t]
	\centering
	\begin{subfigure}{0.49\linewidth}
		\includegraphics[width = \linewidth]{images/Data/stdln_plotsData.pdf}
	\end{subfigure}
	\hfill
	\begin{subfigure}{0.49\linewidth}
		\includegraphics[width = \linewidth]{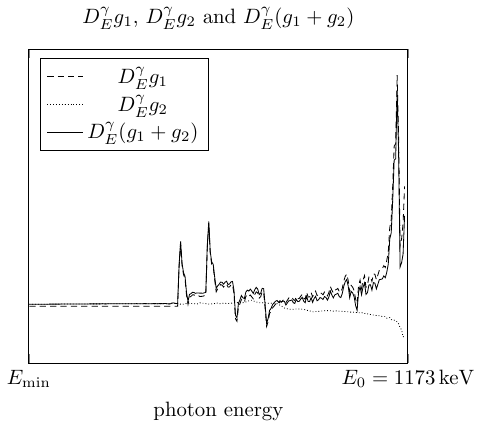}
	\end{subfigure}
	\caption{For the thorax phantom and $E_0 = 1.173$ MeV, we depict the data of a single source-detector pair to show how the differential operator influences the data. Left: unprocessed data, right: the data after applying $D_E^\gamma$.}
	\label{fig:plotsOfDataAndDerivative}
\end{figure}

\Cref{fig:plotsOfDataAndDerivative} shows the result of applying $\mathcal{D}_E^\gamma$ to a measured spectrum. For the computation of the data, the thorax phantom from \cref{sec:modelg1} is used and the measured energy spectrum for one fixed source-detector pair is extracted.
It can be seen that the overall contribution of $g_2$ is fairly large, which explains why reconstructions with forward operator $\mathcal{T}_1$ fail. At the same time, $g_2$ is more spread out on the energy range below the largest value $E_0$ and (in qualitative comparison with $g_1$) it doesn't include any large first or second-order variations. Hence, after application of the derivative in the energy, the main part of the differentiated spectrum belongs to $\mathcal{D}^\gamma_E g_1(\cdot, \mathbf{s},\mathbf{d})$. It is visible how larger variations in $g_1$ (corresponding to contrast jumps in the object) lead to peaks in $\mathcal{D}^\gamma_E g_1(\cdot, \mathbf{s},\mathbf{d})$, which are preserved in the differentiated spectrum, the quantity that can be measured and computed.

We finish the section by noting that both problems \cref{eqn:CSTvariationalForm1} and \cref{eqn:CSTvariationalForm2} admit a solution. Clearly the objective functionals are both convex and bounded from below by 0, so all we need are coercivity results for both, which follow from the structure of the operator $\mathcal{L}_1^{\tilde{n}_e}$.
\begin{corollary}\label{cor:CSTrecProblemExistenceSolution}
	Assume that the object is supported in $\Omega \subset B_1(0)$. Then the problems \cref{eqn:CSTvariationalForm1} and \cref{eqn:CSTvariationalForm2} admit a solution.
\end{corollary}

\begin{proof}
	The TV penalty functional is absolutely one-homogeneous and hence coercive on the complement of its kernel, and the only subspace of $L^2(\Omega)$ on which it vanishes is the set of constant functions $\mathcal{M} = \{f \in L^2(\Omega), f = const. \ \text{a.e.} \}$. Hence coercivity of the objective is ensured if the fidelity term is coercive on $\mathcal{M}$. It suffices to check that $\mathcal{L}_1^{\tilde{n}_e} \chi_\Omega$ resp. $\mathcal{D}_E \mathcal{L}_1^{\tilde{n}_e} \chi_\Omega$ do not vanish, where $\chi_\Omega$ is the indicator function of $\Omega$, i.e. 1 on $\Omega$ and 0 otherwise. We have
	\begin{align*}
		\mathcal{L}_1^{\tilde{n}_e} \chi_\Omega (p, \mathbf{s}, \mathbf{d}) &= \int_{\Omega} \mathcal{W}_1(\tilde{n}_e)(\mathbf{x}, \mathbf{d}, \mathbf{s}) \delta(p - \phi(\mathbf{x}-\mathbf{s}, \mathbf{d}-\mathbf{s})) \mathrm{d} \mathbf{x}\\
		&= \int_{\Omega \cap \mathfrak{C}(\cot^{-1}(p),\mathbf{s},\mathbf{d})} \underset{> 0}{\underbrace{\mathcal{W}_1(\tilde{n}_e)(\mathbf{x}, \mathbf{d}, \mathbf{s})}} \mathrm{d} \mathbf{x} > 0
	\end{align*}
	for any of the circular arcs $\mathfrak{C}(\omega,\mathbf{s},\mathbf{d})$, $p = \cot \omega$ intersecting $\Omega$ (we may assume $\Omega$ has nonzero Lebesgue measure to ensure the integration set is not a null set).\\
	For the version with a derivative in the energy applied to the forward model, we need to show $\mathcal{D}_E \mathcal{L}_1^{\tilde{n}_e} \chi_\Omega \neq 0$. Assume the opposite. Then it holds
	\begin{align*}
		\mathcal{L}_1^{\tilde{n}_e} \chi_\Omega (\cdot, \mathbf{s}, \mathbf{d}) =
		\int_{\Omega \cap \mathfrak{C}(\cot^{-1}(\cdot),\mathbf{s},\mathbf{d})} \mathcal{W}_1(\tilde{n}_e)(\mathbf{x}, \mathbf{d}, \mathbf{s}) \mathrm{d} \mathbf{x} \equiv c(\mathbf{s}, \mathbf{d})
	\end{align*}
	for some constant $c(\mathbf{s},\mathbf{d})$. 
	For any $\mathbf{d} = \mathbf{d}(\theta)$, the circular arc corresponding to $\mathfrak{C}(\cot^{-1}(p),\mathbf{s},\mathbf{d})$ where $\cot^{-1}(p) = \omega = \pi/2 - \theta$ is a subset of $S^1$, the boundary of $B_1(0)$, and thus does not intersect $\Omega$. Therefore, $c(\mathbf{s},\mathbf{d}) = 0$ for every $(\mathbf{s},\mathbf{d})$. But this is a contradiction to the first part of the proof where we showed $\mathcal{L}_1^{\tilde{n}_e} \chi_\Omega \neq 0$.
\end{proof}

Eventually, note that while we only added $g_2$ to the spectrum here, the idea extends to the higher-order terms. According to the concluding remarks about $g_i$ with $i\ge 3$ at the end of \cref{subsec:MappingProperties}, these terms should be smaller in magnitude and at least as smooth as $g_2$. If we succeed in decreasing the impact of $g_2$ sufficiently by differentiating $g_1 + g_2$, we can reasonably expect the procedure to work on the complete spectrum.

\subsection{Extension to Polychromatic Sources}\label{subsec:Polychromatic}
For simplicity, until here we assumed the source to be monochromatic. Consider now the case where the source emits photons at $K$ different energy levels with total initial intensity $I$ and a finite, '{}small'{} number $K$. Then $I = I_1 + \dots + I_K$, where $I_k$ is the total intensity of all photons of initial energy $E_0^k$, $k = 1,\dots, K$. 

Recall that our modeling of the first-order scattering emerged from the representation of the number of photons $N$ reaching a detector at $\mathbf{d}$ after being scattered at $\mathrm{d} \mathbf{x}$ through
\begin{equation*}
	\mathrm{d}^2 N(\mathbf{s},\mathbf{d},\mathbf{x}) = I r_e^2 A_{E_0}(\mathbf{s},\mathbf{x}) A_{E_\omega}(\mathbf{x}, \mathbf{d}) P(\omega, E_0) n_e(\mathbf{x})\, \mathrm{d} \mathbf{x}\,\mathrm{d}\mathbf{s}
\end{equation*}
for a monochromatic source. Note that we write $P(\omega, E_0)$ instead of $P(\omega)$ indicating the dependence of the Klein-Nishina differential cross section on the photon energy. Decomposing the effect for a polychromatic source with finitely many initial energies into the single energy levels, this model remains accurate if we consider the $K$ levels separately, i.e.
\begin{equation*}
	\mathrm{d}^2 N(\mathbf{s},\mathbf{d},\mathbf{x}) = I r_e^2 \left(\sum_{k=1}^{K} c_k P(\omega, E_0^k) A_{E_0^k}(\mathbf{s},\mathbf{x}) A_{E_\omega^k}(\mathbf{x}, \mathbf{d})\right) n_e(\mathbf{x})\, \mathrm{d} \mathbf{x}\,\mathrm{d}\mathbf{s},
\end{equation*}
where $c_k := I_k/I$ is a weight factor. As the integration over the possible scattering sites is linear, the first-order scattered spectrum for a polychromatic source is a weighted sum of operators $\mathcal{T}_1^{E_0^k}$ for every initial energy level $E_0^k$. Similarly, given a prior electron density map $\tilde{n}_e$, we can define a polychromatic version of the linearized operator by
\begin{equation*}
	\mathcal{T}_1^{\mathrm{pol}} = \sum_{k=1}^{K} c_k \mathcal{T}_1^{E_0^k}\quad\textrm{and}\quad
	\mathcal{L}_1^{\tilde{n}_e, \mathrm{pol}} = \sum_{k=1}^{K} c_k \mathcal{L}_1^{\tilde{n}_e, E_0^k}.
\end{equation*}
The mapping properties proved in \cref{subsec:MappingProperties} hold for the operators $\mathcal{L}_1^{\tilde{n}_e, E_0^k}$ and therefore extend to their weighted sum $\mathcal{L}_1^{\tilde{n}_e, \mathrm{pol}}$. Equivalently, the modeling can be extended for the second-order scattered data, allowing us to define $\mathcal{T}_2^{\mathrm{pol}}$ and $\mathcal{L}_2^{\tilde{n}_e, \mathrm{pol}}$ as weighted sums and prove suitable mapping properties. The reconstruction task \cref{eqn:CSTvariationalForm2} remains the same, with $\mathcal{L}_1^{\tilde{n}_e}$ replaced by $\mathcal{L}_1^{\tilde{n}_e, \mathrm{pol}}$. We describe a possible scanner setup and show the resulting density reconstructions in the next section.

\section{Density Reconstructions for a Fan-Beam Geometry}\label{sec:Results}
We provide results for the setting of a monochromatic source first. In practice, potential sources emit photons at several energy levels. We discuss what properties the source would need to satisfy for our algorithm to work and then depict results in the case of a Cobalt-60 gamma ray source.
\subsection{Monochromatic Source}\label{subsec:resultsMonochromatic}
Following the discussions in \cref{subsec:Algorithm}, we propose to combine CST and CT in the same scanning geometry. As we wish to acquire several viewing angles (resp. positions of the source), a standard fan beam setup appears sensible. The major difference to a classical CT setup is that we assume the radiation detectors to be energy resolving and non-collimated so that our models of $g_1$ and $g_2$ are applicable.

\begin{figure}
	\centering
	\includegraphics[width=0.4\textwidth]{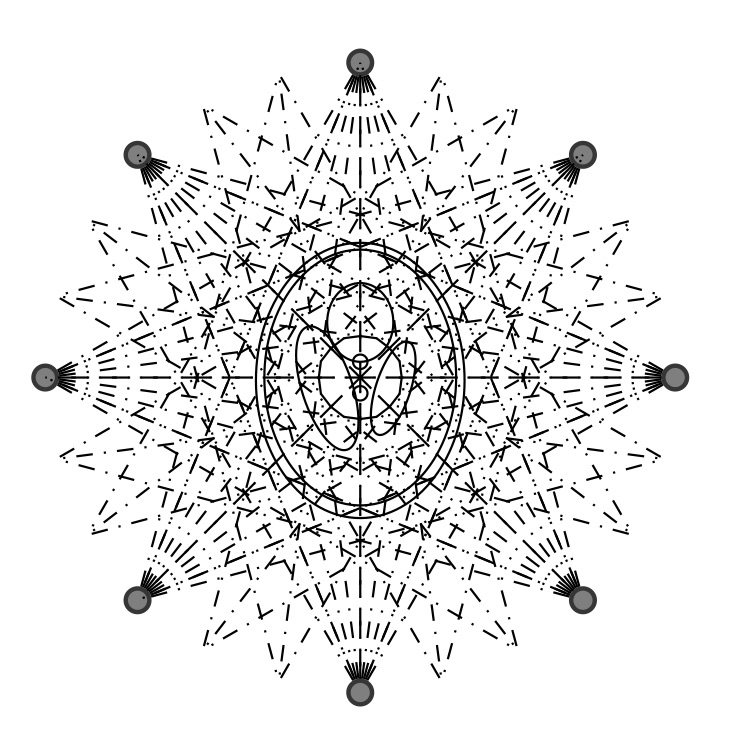}
	\caption{Scanning of the object for 8 angular views.}
	\label{fig:ScanGeometry}
\end{figure}

As was mentioned in the description of the reconstruction algorithm, the procedure is carried out with a very small number of source and detector positions. This is enabled by the energy as additional dimension of information in the CST step. Assume a setup with 16 angular views, see \Cref{fig:ScanGeometry} for 8 angular views. For every angle, the radiation behind the object is measured by 32 detectors distributed equidistantly around a half circle opposite to the source. Particularly, note that this setup is an instance of the geometry assumed in the Sobolev mapping results \Cref{thm:sobolevtheoreml1,thm:sobolevtheoreml2}. Regarding the radiation source, we assume a total emission intensity of $10^7$ photons per viewing angle and an initial photon energy level of $1.173\,$MeV. The detectors measure radiation at 256 energy levels in the range $0.355\,$MeV to $1.173\,$MeV, leading to a necessary energy resolution $\Delta E \approx 3.2\,$keV. The nature of the source is incorporated in the simulation of the data by computing Poisson distributed values $g(E,\mathbf{s},\mathbf{d}) \sim \Pois((\mathcal{T}_1 + \mathcal{T}_2)(n_e^{\mathrm{ex}})(E,\mathbf{s},\mathbf{d}))$.

The final reconstruction procedure consists in first solving the sparse data CT step \cref{eqn:CTLimitedDataMinimization} and then, using the prior reconstruction, solving \cref{eqn:CSTvariationalForm2}. The problem is converted into the discrete setting by decomposing the image $f$ in a suited pixel basis giving a vector representation $\bar{f} \in \mathbb{R}^n$ and computing a sparse matrix representation $L_1^{\tilde{n}_e} \in \mathbb{R}^{m \times n}$ of $\mathcal{L}_1^{\tilde{n}_e}$ with respect to the basis functions.
The gradient in $\TV$ is substituted with forward finite differences. As in \cite{Acar94, Almansa08}, we obtain a regularized isotropic TV of the form
\begin{equation*}
	J_\beta(\bar{f}) := \sum_{1 \le i,j \le N} \sqrt{\abs{\nabla \bar{f}(i,j)}^2 + \beta},
\end{equation*}
where commonly a small number $\beta > 0$ is added to the gradient to avoid computational instabilities in constant regions.

The computationally intense part of the reconstruction is the generation of the discrete representation of the linearized forward operator $L_1^{\tilde{n}_e}$. Every entry of the $N^2 \times N_{\textrm{det}}N_{\textrm{sou}}N_{\textrm{E}}$ matrix requires the evaluation of the integral (11) along a circular arc, which is done by interpolating the values of the weight function $\mathcal{W}_1(\tilde{n}_e)$ at quadrature points $\mathbf{x}_i$ inside a basis function support. We assume $N_{\textrm{det}}N_{\textrm{sou}}N_{\textrm{E}} = \mathcal{O}(N^2)$ which gives a total amount of $\mathcal{O}(N^4)$ operations. If we were to include $L_2^{\tilde{n}_e}$ into the forward model, then the complexity would increase to $\mathcal{O}(N^6)$, which further motivates the solution of the problem using forward models including only sinngle scattering. In 3D, the costs of computing $L_1^{\tilde{n}_e}$ and $L_2^{\tilde{n}_e}$ are $\mathcal{O}(N^6)$ and $\mathcal{O}(N^9)$ respectively.

\begin{figure}[t]
	\centering
	\begin{subfigure}{0.31\textwidth}
		\centering
		\includegraphics[width=\textwidth]{images/Thorax/groundTruthThorax.png}
		\subcaption{}
		\label{subfig:recMCgroundTruthThorax}
	\end{subfigure}
	\hspace{0.1cm}
	\begin{subfigure}{0.31\textwidth}
		\centering
		\includegraphics[width=\textwidth]{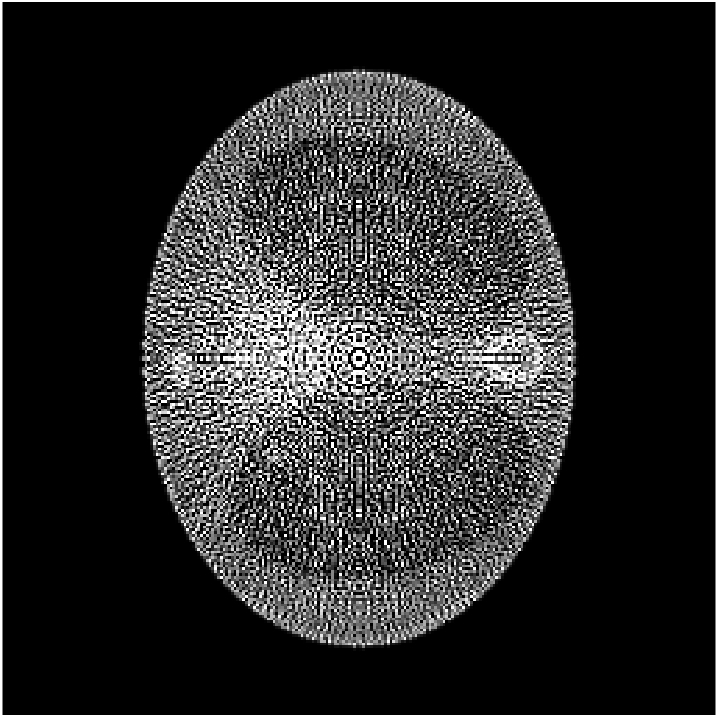}
		\subcaption{}
		\label{subfig:recMCThoraxnoDnoTV}
	\end{subfigure}
	\hspace{0.1cm}
	\begin{subfigure}{0.31\textwidth}
		\centering
		\includegraphics[width=\textwidth]{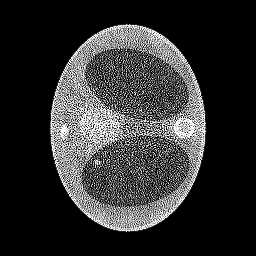}
		\subcaption{}
		\label{subfig:recMCThoraxDnoTV}
	\end{subfigure}\\
	\begin{subfigure}{0.31\textwidth}
		\centering
		\includegraphics[width=\textwidth]{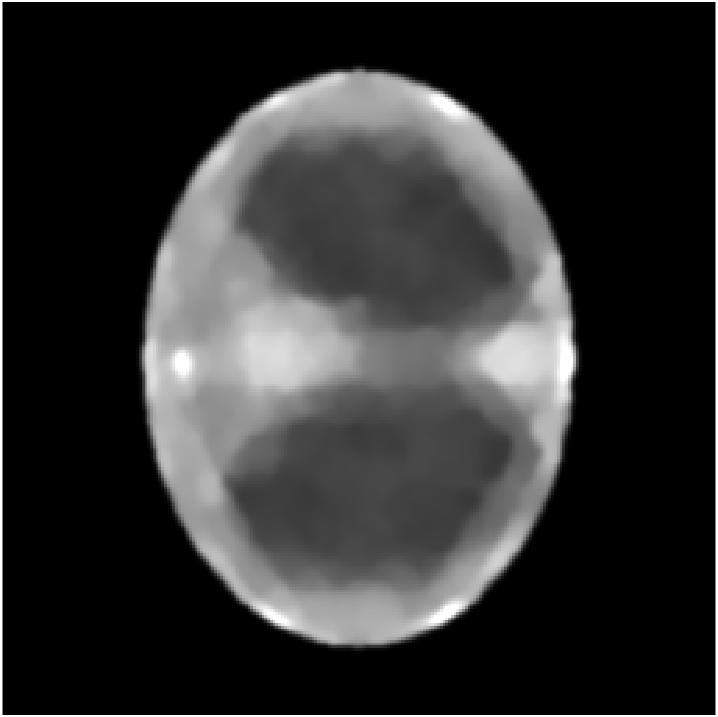}
		\subcaption{}
		\label{subfig:recMCThoraxCTprior}
	\end{subfigure}
	\hspace{0.1cm}
	\begin{subfigure}{0.31\textwidth}
		\centering
		\includegraphics[width=\textwidth]{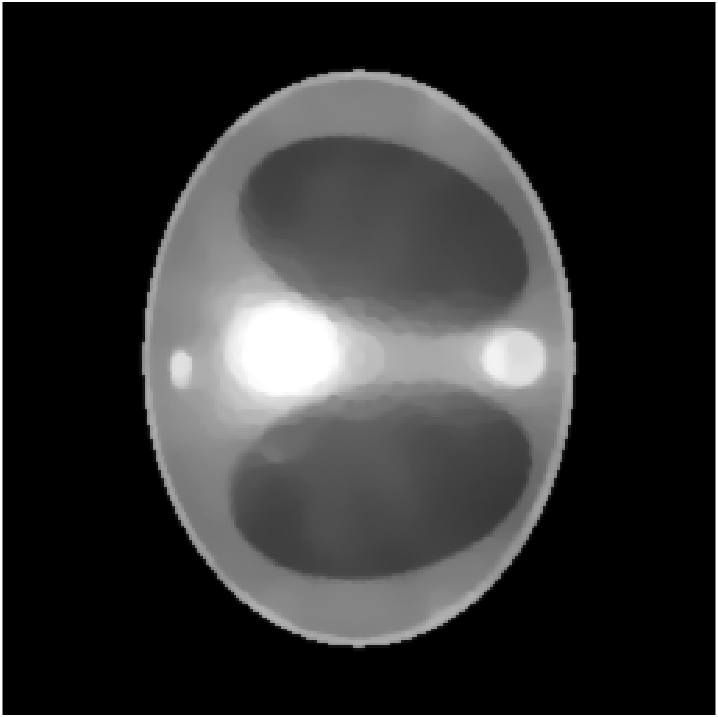}
		\subcaption{}
		\label{subfig:recMCThoraxnoDTV}
	\end{subfigure}
	\hspace{0.1cm}
	\begin{subfigure}{0.31\textwidth}
		\centering
		\includegraphics[width=\textwidth]{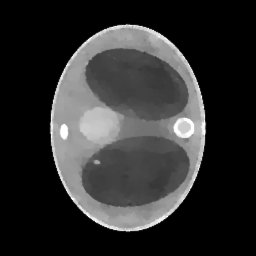}
		\subcaption{}
		\label{subfig:recMCThoraxDTV}
	\end{subfigure}
	\caption{(\ref{subfig:recMCgroundTruthThorax}): Ground truth phantom. (\ref{subfig:recMCThoraxCTprior}): prior reconstruction from sparse data CT step \cref{eqn:CTLimitedDataMinimization}. (\ref{subfig:recMCThoraxnoDnoTV}) and (\ref{subfig:recMCThoraxnoDTV}): Solutions of problem \cref{eqn:CSTvariationalForm1}, $\lambda = 0$ in (\ref{subfig:recMCThoraxnoDnoTV}) and $\lambda > 0$ tuned in (\ref{subfig:recMCThoraxnoDTV}).
		(\ref{subfig:recMCThoraxDnoTV}) and (\ref{subfig:recMCThoraxDTV}): Solutions of problem \cref{eqn:CSTvariationalForm2} with TV parameters $\lambda = 0$ (least-squares fit) in (\ref{subfig:recMCThoraxDnoTV}) and $\lambda > 0$ tuned by the L-curve method in (\ref{subfig:recMCThoraxDTV}).}
	\label{fig:ReconstructionsThoraxMC}
\end{figure}

The method is tested for two synthetic phantoms; firstly, the thorax phantom already used in \cref{sec:modelg1}.
\Cref{fig:ReconstructionsThoraxMC} shows how the algorithm succeeds in decreasing the impact of the higher-order scattering. The ground truth is displayed in \cref{subfig:recMCgroundTruthThorax}. Expectably, the CT reconstruction $\tilde{n}_e$ (\cref{subfig:recMCThoraxCTprior}) is not accurate enough, but can be used as a prior to estimate the nonlinear weight function $\mathcal{W}_1(\tilde{n}_e)$. Using $\tilde{n}_e$, the CST reconstruction is computed. For comparison, we give both the solution of \cref{eqn:CSTvariationalForm1} (without $\mathcal{D}_E$) and \cref{eqn:CSTvariationalForm2} (with $\mathcal{D}_E$). As in \cref{sec:modelg1}, \cref{eqn:CSTvariationalForm1} cannot yield a useful reconstruction. The minimizer when $\mathrm{R}$ is the Kullback-Leibler divergence \eqref{eq:kullbackLeiblerDivDef} and $\lambda = 0$ is very noisy (\cref{subfig:recMCThoraxnoDnoTV}) and using TV regularization, some noise can be filtered out, but only at the cost of losing small details of lower contrast (\cref{subfig:recMCThoraxnoDTV}). As desired, applying $\mathcal{D}_E$ to the data reduces the noise level, see \cref{subfig:recMCThoraxDnoTV} and \cref{subfig:recMCThoraxDTV} (with TV regularization). After applying $\mathcal{D}_E$, we use the $L_2$ norm as data fidelity measure. Densities and contrasts are accurately recovered and previously vanished details can be correctly located.

\begin{figure}[t]
	\centering
	\begin{subfigure}{0.31\textwidth}
		\centering
		\includegraphics[width=\textwidth]{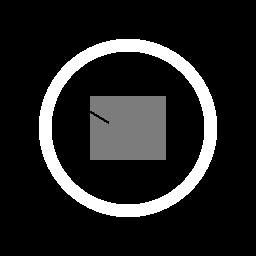}
		\subcaption{}
		\label{subfig:recMCgroundTruthAlu}
	\end{subfigure}
	\hspace{0.1cm}
	\begin{subfigure}{0.31\textwidth}
		\centering
		\includegraphics[width=\textwidth]{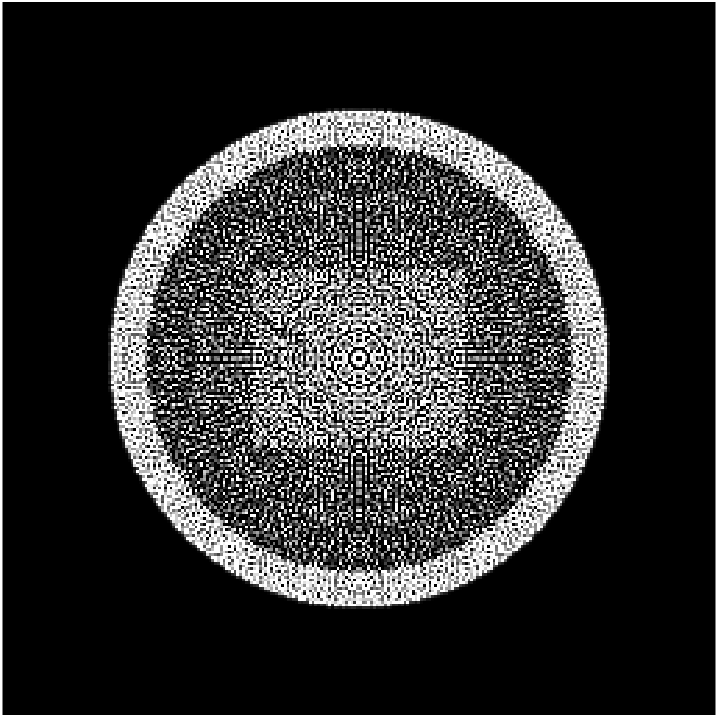}
		\subcaption{}
		\label{subfig:recMCAlunoDnoTV}
	\end{subfigure}
	\hspace{0.1cm}
	\begin{subfigure}{0.31\textwidth}
		\centering
		\includegraphics[width=\textwidth]{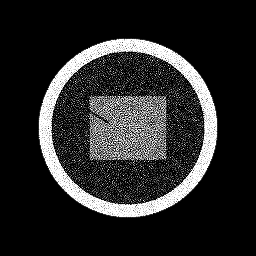}
		\subcaption{}
		\label{subfig:recMCAluDnoTV}
	\end{subfigure}
	\\
	\begin{subfigure}{0.31\textwidth}
		\centering
		\includegraphics[width=\textwidth]{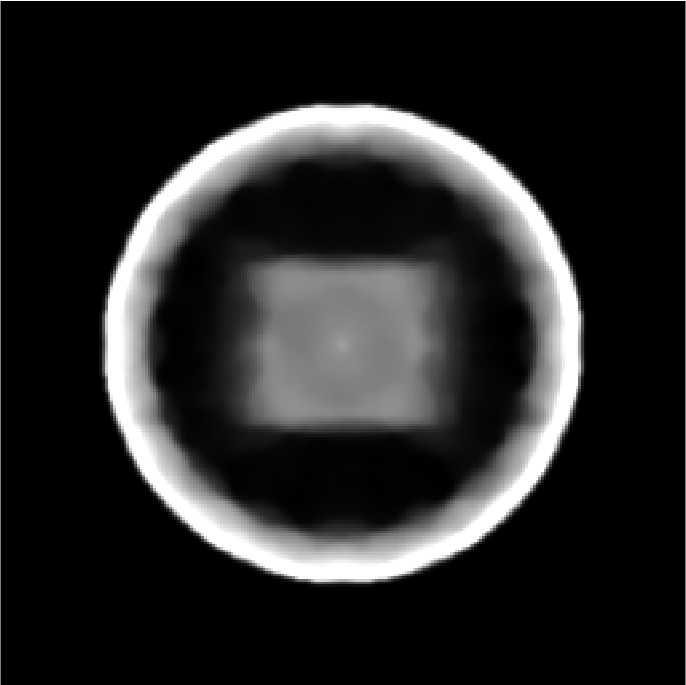}
		\subcaption{}
		\label{subfig:recMCAluCTprior}
	\end{subfigure}
	\hspace{0.1cm}
	\begin{subfigure}{0.31\textwidth}
		\centering
		\includegraphics[width=\textwidth]{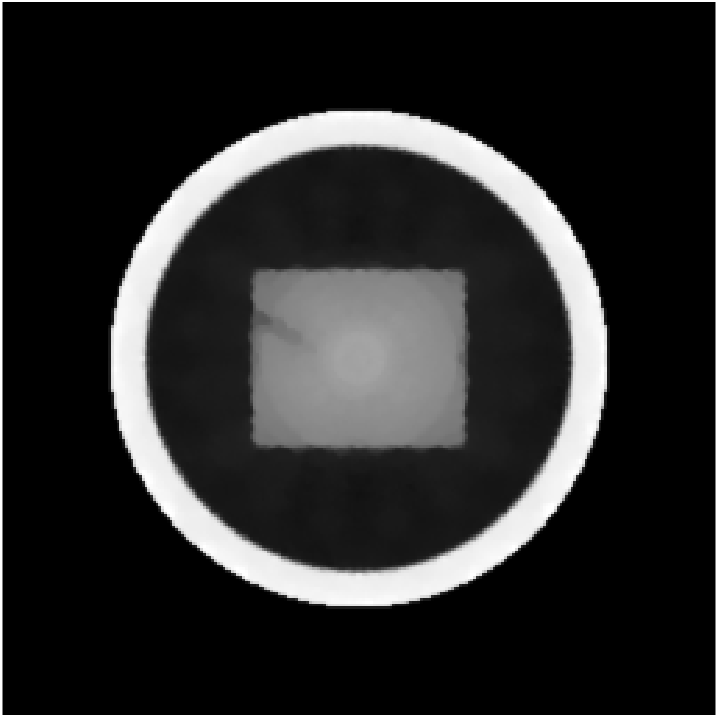}
		\subcaption{}
		\label{subfig:recMCAlunoDTV}
	\end{subfigure}
	\hspace{0.1cm}
	\begin{subfigure}{0.31\textwidth}
		\centering
		\includegraphics[width=\textwidth]{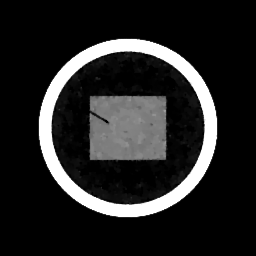}
		\subcaption{}
		\label{subfig:recMCAluDTV}
	\end{subfigure}
	\caption{(\ref{subfig:recMCgroundTruthAlu}): The second phantom with an aluminium ring. (\ref{subfig:recMCAluCTprior}): prior reconstruction from sparse data CT step \cref{eqn:CTLimitedDataMinimization}. (\ref{subfig:recMCAlunoDnoTV}) and (\ref{subfig:recMCAlunoDTV}): Solutions of problem \cref{eqn:CSTvariationalForm1} with $\lambda = 0$ (\ref{subfig:recMCAlunoDnoTV}) and $\lambda > 0$ tuned (\ref{subfig:recMCAlunoDTV}).
		(\ref{subfig:recMCAluDnoTV}) and (\ref{subfig:recMCAluDTV}): Solutions of problem \cref{eqn:CSTvariationalForm2} with $\lambda = 0$ (least-squares fit) in (\ref{subfig:recMCAluDnoTV}) and $\lambda > 0$ tuned in (\ref{subfig:recMCAluDTV}).}
	\label{fig:Al_recs_MC}
\end{figure}

\begin{figure}[t]
	\centering
	\begin{subfigure}{0.31\textwidth}
		\centering
		\includegraphics[width=\textwidth]{images/Pipe/Iron/MC/groundTruth.png}
		\subcaption{}
		\label{subfig:recMCgroundTruthIron}
	\end{subfigure}
	\hspace{0.1cm}
	\begin{subfigure}{0.31\textwidth}
		\centering
		\includegraphics[width=\textwidth]{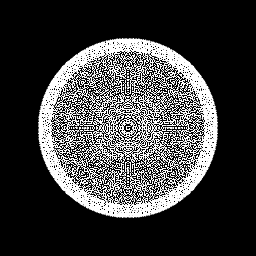}
		\subcaption{}
		\label{subfig:recMCIronnoDnoTV}
	\end{subfigure}
	\hspace{0.1cm}
	\begin{subfigure}{0.31\textwidth}
		\centering
		\includegraphics[width=\textwidth]{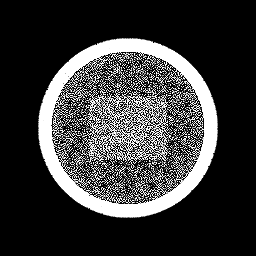}
		\subcaption{}
		\label{subfig:recMCIronDnoTV}
	\end{subfigure}
	\\
	\begin{subfigure}{0.31\textwidth}
		\centering
		\includegraphics[width=\textwidth]{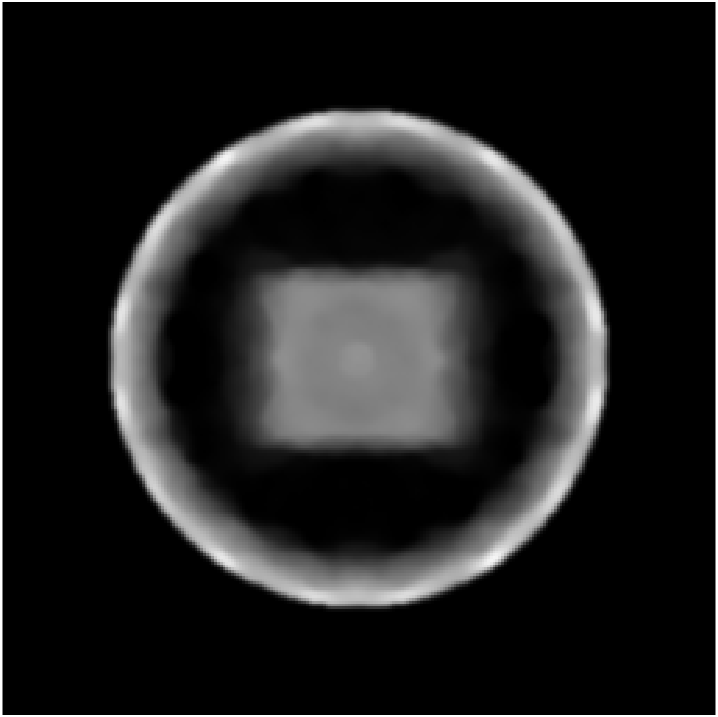}
		\subcaption{}
		\label{subfig:recMCIronCTprior}
	\end{subfigure}
	\hspace{0.1cm}
	\begin{subfigure}{0.31\textwidth}
		\centering
		\includegraphics[width=\textwidth]{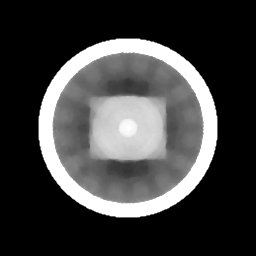}
		\subcaption{}
		\label{subfig:recMCIronnoDTV}
	\end{subfigure}
	\hspace{0.1cm}
	\begin{subfigure}{0.31\textwidth}
		\centering
		\includegraphics[width=\textwidth]{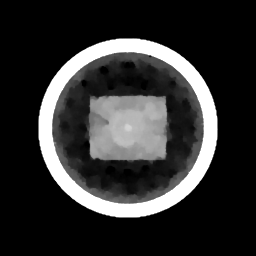}
		\subcaption{}
		\label{subfig:recMCIronDTV}
	\end{subfigure}
	\caption{(\ref{subfig:recMCgroundTruthIron}): The second phantom, now with an iron ring. (\ref{subfig:recMCIronCTprior}): prior reconstruction from sparse data CT step \cref{eqn:CTLimitedDataMinimization}. (\ref{subfig:recMCIronnoDnoTV}) and (\ref{subfig:recMCIronnoDTV}): Solutions of problem \cref{eqn:CSTvariationalForm1} with $\lambda = 0$ in (\ref{subfig:recMCIronnoDnoTV}) and $\lambda > 0$ tuned in (\ref{subfig:recMCIronnoDTV}).
		(\ref{subfig:recMCIronDnoTV}) and (\ref{subfig:recMCIronDTV}): Solutions of problem \cref{eqn:CSTvariationalForm2} with TV parameters $\lambda = 0$ (least-squares fit) in (\ref{subfig:recMCIronDnoTV}) and $\lambda > 0$ tuned in (\ref{subfig:recMCIronDTV}).}
	\label{fig:Iron_recs_MC}
\end{figure}

We further use two variants of a second phantom (\cref{subfig:recMCgroundTruthAlu,subfig:recMCgroundTruthIron}), consisting of materials whose densities are more likely to occur in industrial contexts. Its outside is a metal ring of outer radius 3.5 cm and inner radius 3 cm. Inside the ring is a rectangle of size $1.5 \times 1.25 \cm^2$ from a material of smaller density, in our model the plastic type polyethylene (PE). The plastic is cracked at one point and we wish to correctly recover the object so that the crack of width 1 mm is visible and can be located correctly. To show the limitations of the method and the impact of the electron density (leading to more scattering and larger nonlinear effects) of the scanned material, the outer ring is hereby once assumed to be made from aluminium (moderate density) and once from iron (higher density). 

\Cref{fig:Al_recs_MC} shows the reconstructions in the case of an aluminium ring on a $256\times 256$ pixel image (one pixel corresponds to $0.20 \times 0.20 \mm^2$). In \cref{fig:Iron_recs_MC} the same results but with an iron ring are displayed. From a sparse-view CT step (16 source positions and 32 detectors), the rough contours of the object are reconstructed well if we apply large $\TV$ parameters, but there are artifacts and the crack is clearly invisible, see \cref{subfig:recMCAluCTprior,subfig:recMCIronCTprior}. Also, the iron already creates greater artifacts so that larger $\TV$ parameters have to be chosen to obtain a smooth image.
Taking $\mathcal{L}_1^{\tilde{n}_e}$ as the forward model, the different effects of aluminium and iron on the data become clearer: The higher attenuation of iron leads to more noise and circular artifacts aligned with the source positions, \cref{subfig:recMCIronnoDnoTV} is worse than the aluminium case in \cref{subfig:recMCAlunoDnoTV}. Using TV-regularization, the crack is vaguely visible in the aluminium phantom, but not in the iron version (\cref{subfig:recMCAlunoDTV,subfig:recMCIronnoDTV}). Including the derivative greatly improves the quality of both reconstructions (\cref{subfig:recMCAluDnoTV,subfig:recMCAluDTV,subfig:recMCIronDnoTV,subfig:recMCIronDTV}), although only the aluminium version with tuned TV parameter appears accurate enough to exactly determine the nature of the crack in the plastic.

\subsection{Polychromatic Source}
We adopt the theory from \cref{subsec:Polychromatic} and extend the model to a polychromatic source. Naturally, several requirements on the source material arise. 

Most importantly, it is necessary that we are able to separate $g_0$ from the rest of the measured spectrum. In the case of a monochromatic source, this is easy since $g_0$ belongs to exactly those photons that have energy $E_0$. Assuming sufficient energy resolution and no detector defects (e.g. photons whose energies add up to $\approx E_0$ that hit the detector at the same time), the peak measured at $E_0$ can just be extracted giving $g_0$. In a polychromatic setting, this strategy can still work for the highest initial energy level $E_0^{\mathrm{max}}$. For the initial energy levels below this, the job is more delicate: We can only estimate the component $g_0$ coming from initial energy levels $E_0^k \neq E_0^{\max}$ using, e.g., an interpolation of the spectrum around $E_0^k$ or by the total intensity deposited in the detector at energy $E_0^{\mathrm{max}}$. Effects like beam hardening further complicate this task. Therefore, we wish to keep the number of energy levels as small as possible, where the initial energy levels are sufficiently distinct from each other. From a computational point of view, a small number of energy levels is further encouraged bearing in mind that the computation of weight matrices $L_1^{\tilde{n}_e,E_0^k}$, $k=1,\dots,K$ is cumbersome.

We assume a radiation source consisting of Cobalt-60. In a $\beta^-$ decay, the $\textrm{Co}^{60}_{27}$ nuclei decay into $\textrm{Ni}_{28}^{60}$, causing the latter, excited nuclei to emit photons of only two energy levels with the values $E_0^1 = 1.173$\,MeV, $E_0^2 = 1.332$\,MeV.
60Co was a popular material in radiotherapy due to its relatively long half life and with its only two initial energies that are sufficiently distinct it fits the requirements we have put on a potential source material.

\begin{figure}[t]
	\centering
	\begin{subfigure}{0.31\textwidth}
		\centering
		\includegraphics[width=\textwidth]{images/Thorax/groundTruthThorax.png}
		\subcaption{}
		\label{subfig:recPCgroundTruthThorax}
	\end{subfigure}
	\hfill
	\begin{subfigure}{0.31\textwidth}
		\centering
		\includegraphics[width=\textwidth]{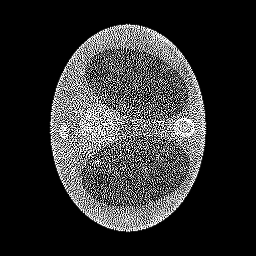}
		\subcaption{}
		\label{subfig:recPCThoraxDnoTV}
	\end{subfigure}
	\hfill
	\begin{subfigure}{0.31\textwidth}
		\centering
		\includegraphics[width=\textwidth]{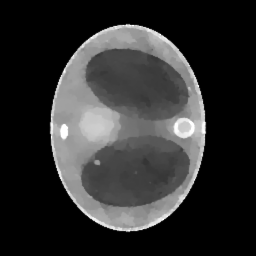}
		\subcaption{}
		\label{subfig:recPCThoraxDTV}
	\end{subfigure}
	\caption{(\ref{subfig:recPCgroundTruthThorax}): Ground truth phantom. (\ref{subfig:recPCThoraxDnoTV}) and (\ref{subfig:recPCThoraxDTV}): Solution of \cref{eqn:CSTvariationalForm2} with $\lambda = 0$ in (\ref{subfig:recPCThoraxDnoTV}) (least-squares fit) and $\lambda > 0$ tuned in (\ref{subfig:recPCThoraxDTV}).}
	\label{fig:ReconstructionsThoraxPC}
\end{figure}

\begin{figure}[t]
	\centering
	\begin{subfigure}{0.31\textwidth}
		\centering
		\includegraphics[width=\textwidth]{images/Pipe/Iron/MC/groundTruth.png}
		\subcaption{}
		\label{subfig:recPCgroundTruthAlu}
	\end{subfigure}
	\hfill
	\begin{subfigure}{0.31\textwidth}
		\centering
		\includegraphics[width=\textwidth]{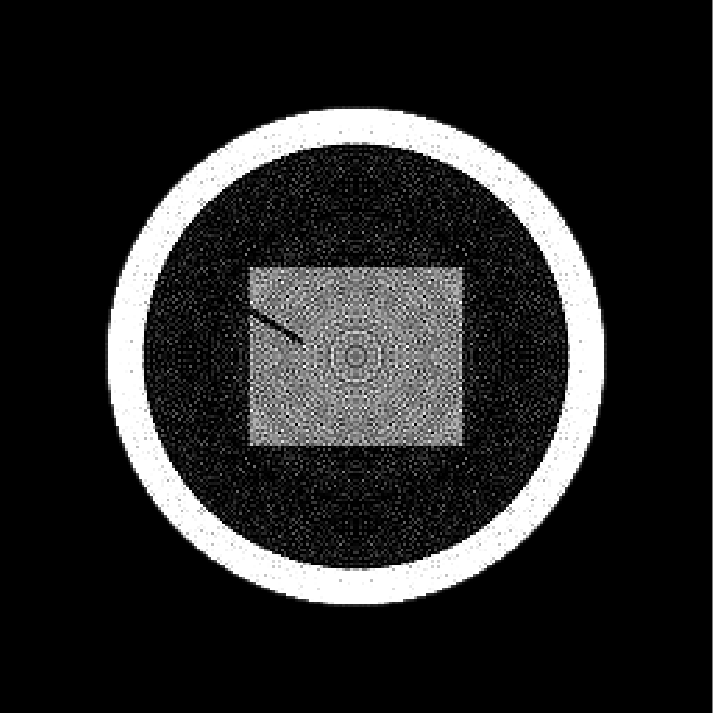}
		\subcaption{}
		\label{subfig:recPCAluDnoTV}
	\end{subfigure}
	\hfill
	\begin{subfigure}{0.31\textwidth}
		\centering
		\includegraphics[width=\textwidth]{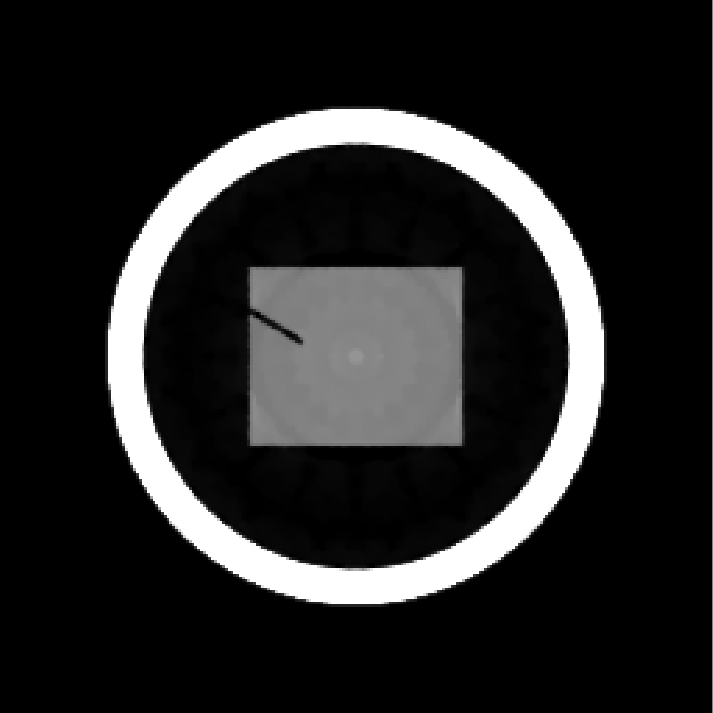}
		\subcaption{}
		\label{subfig:recPCAluDTV}
	\end{subfigure}
	\caption{Results for the aluminium object (\ref{subfig:recPCgroundTruthAlu}): Ground truth phantom. (\ref{subfig:recPCAluDnoTV}) and (\ref{subfig:recPCAluDTV}): Solution of \cref{eqn:CSTvariationalForm2} with $\lambda = 0$ in (\ref{subfig:recPCAluDnoTV}) (least-squares fit) and $\lambda > 0$ tuned in (\ref{subfig:recPCAluDTV}).}
	\label{fig:ReconstructionsAluPC}
\end{figure}

Apart from the computation of a second weight matrix, the reconstructions are computed the same way as in the monochromatic case. We only give the solutions of \cref{eqn:CSTvariationalForm2} for the thorax phantom (\cref{fig:ReconstructionsThoraxPC}) and the second phantom in the aluminium case (\cref{fig:ReconstructionsAluPC}).
The results are of equally convincing quality as in the case of a monochromatic source. With prior from the initial sparse-view CT step, the forward model is sufficiently accurate to recover both phantoms with few artifacts from the data $g_1$. Using a differentiation step, the impact of $g_2$ and possibly higher-order terms on the reconstruction are decreased due to their smoothness. TV regularization is used in order to take care of the remaining noise generated by the Poisson process.

\section{Conclusion \& Discussion}\label{sec:conclusion}
In this paper we proposed a reconstruction strategy for two-dimensional tomographic slices based on a modality which combines computerized tomography and Compton scattering tomography. We considered a standard fan-beam geometry and assumed the detectors are able to measure the energy spectrum of the incoming photons. Due to the spatial sparsity of the overall architecture, the reconstruction of the attenuation map using CT-scan data is of poor quality but provides still a prior estimate of the electron density, which is used to refine the modeling used in CST. Focusing on first-order and second-order scattering, we derived and studied operators describing the polychromatic spectrum. Linearizations of these operators were studied as examples of the general class of Fourier integral operators in order to prove suitable mapping properties.

Due to the complexity of the terms of higher-order in the spectrum, we first restricted the inverse problem to the modeling of the first-order scattering, but this simple approach turned out to be unsatisfactory as the higher-order scattered data bring a large distortion of the spectrum. However, the mapping properties reveal that the contribution of multiple scattering is smoother in comparison with the first-order scattering. We thus exploited the structure of the spectrum by applying a differential operator to the data. The resulting algorithm was validated for phantoms characterizing potential applications in medicine and non-destructive testing. The results demonstrate how a qualitative analysis of terms in the spectrum can lead to improved images and solve the problem of multiple scattering in CST.

Moreover, we extended the model to polychromatic sources of few, sufficiently distinct emitted energy levels. Due to the flexibility of the algebraic approach, this can be done by a simple weighted sum of the operators of distinct energy levels of a radioactive gamma ray source like Cobalt-60, requiring hardly any additional effort over the monochromatic case, but larger computational time.

Just as in algebraic reconstruction methods in CT, for images of appropriate resolution the  projection matrix is large. Hence, computational time must be taken into account as a relevant factor in real-time applications. Our discussion of physical limitations of the method revealed two further bottlenecks in the method. The quality of the prior obtained from the initial CT step is crucial as it influences the consistency between $\mathcal{T}_1$ and its linearization $\mathcal{L}_1^{\tilde{n}_e}$. Further, the comparison of iron and aluminium phantoms showed how a larger density also effects the total numbers of once or twice scattered photons and hence the reconstruction quality. The size of a scanned object and the densities of materials in it therefore have a major influence on the image quality. 

Future work will be dedicated to apply the method to e.g. Monte-Carlo simulated or real data from a scanner to validate the model. Further, the model will be improved by incorporating the physical size of source and detectors in order to approach a realistic setting. In the same spirit, the algorithm can be extended to the 3D case, as Compton scattering is naturally a three-dimensional process. The extension to 3D appears promising though computationally more difficult, since similar smoothness properties of the spectral data were proven in \cite{Rigaud21}.

\section{\texorpdfstring{Proofs - Modelling of $\mathcal{T}_2$}{Proofs - Modelling of T2}}\label{pro:appxModelT2}
\begin{proof}[Proof of Lemma \ref{lem:SecondScatteringPoint}]
	After the photon is scattered at $\mathbf{x}$, its trajectory is a straight line that forms an angle of $\pi-\omega_1$ with the incoming beam direction $\mathbf{x} - \mathbf{s}$. Thus, its trajectory is part of the cone opening at $\mathbf{x}$ with direction $\mathbf{x} - \mathbf{s}$ and aperture $\omega_1$. As we only consider the two-dimensional case, the cone simplifies to two possible directions, resembling the shape of a "{}V"{}, see the dotted lines in figure \ref{fig:SecondOrderScGeometry}. We simplify notation by introducing a convention of a signed angle $\omega_1$, restricting $\omega_1$ not to $[0, \pi]$ but to $[-\pi, \pi]$, although both $\pm \omega_1$ yield the same energy $E_1$ after the first scatter event. The photon's trajectory after the first scattering is then a subset of the line
	\begin{equation}\label{eqn:coneParametrization}
		\mathfrak{L}(\omega_1, \mathbf{x}, \mathbf{s}) = \left\{ \mathbf{x} + t \begin{pmatrix} \sin (\omega_1 - \beta + \pi/2)\\ \cos(\omega_1 - \beta + \pi/2)\end{pmatrix},\ t\ge 0 \right\}
	\end{equation}
	with $\omega_1 \in [-\pi, \pi]$, where $\beta \in \left[0,2\pi\right)$ is defined by $\frac{\mathbf{x}-\mathbf{s}}{\norm{\mathbf{x}-\mathbf{s}}} = (\cos \beta, \sin \beta)^T$.
	
	Fixing the energy $E$ at $\mathbf{d}$, as in \cref{eqn:sumCosines} we define $\omega_2(\omega_1) = \arccos(\lambda(E) - \cos(\omega_1))$. Interpreting $\mathbf{x}$ as a radiation source, we know, equivalently as in the modeling of $g_1$, that $\mathbf{x}$ is a point on one of two symmetric circular arcs that can be parametrized by
	\begin{align}
		\mathfrak{C}(\omega_2, \mathbf{d}, \mathbf{x}) &= \left\{ \mathbf{x} + \norm{\mathbf{d}-\mathbf{x}} \frac{\sin(\omega_2 - \alpha)}{\sin \omega_2} R \begin{pmatrix}
			\pm \sin \alpha\\ \cos \alpha
		\end{pmatrix},\ \alpha \in [0, \omega_2] \right\}\notag\\
		&=: \left\{ \mathbf{x} + r(\omega_2, \alpha) d^{\pm}(\alpha),\ \alpha \in [0, \omega_2] \right\}\label{eqn:torusParametrization}.
	\end{align}
	where $R$ is an orthogonal rotation matrix mapping $(0,1)^T$ to $(\mathbf{d} - \mathbf{x})/\norm{\mathbf{d} - \mathbf{x}}$ and the $\pm$-sign distinguishes the two separated circular arcs. The second formulation includes a unit vector $d^{\pm}(\alpha) := R (\pm \sin \alpha, \cos \alpha)^T $ and a corresponding scaling factor $r(\omega_2, \alpha) := \norm{\mathbf{d}-\mathbf{x}} \sin(\omega_2 - \alpha)/\sin \omega_2$.\\
	We can now compute the possible intersection points $\mathbf{y}$ by comparing the parametrizations \eqref{eqn:coneParametrization} and \eqref{eqn:torusParametrization}.
	Subtracting $\mathbf{x}$ and applying $R^{-1}$, we obtain that 
	\begin{equation}\label{eqn:comparisonParametrizations}
		t\eta := t R^{-1} \begin{pmatrix} \sin (\omega_1 - \beta + \pi/2)\\ \cos(\omega_1 - \beta + \pi/2)\end{pmatrix} = r(\omega_2, \alpha) \begin{pmatrix}
			\pm \sin \alpha\\ \cos \alpha
		\end{pmatrix}
	\end{equation}
	at an intersection point, hereby $\eta$ only depends on the location of $\mathbf{x}$ and $\mathbf{d}$. Further, it is well-defined and can be computed easily as $R^{-1} = R^T$, the second component is given by
	\begin{equation}\label{eqn:etaTwo}
		\cos \alpha = \eta_2 = R_{1,2} \sin (\omega_1 - \beta + \pi/2) + R_{2,2} \cos(\omega_1 - \beta + \pi/2).
	\end{equation}
	
	Using trigonometric identities, we can rewrite $r$ only depending on the value of $\omega_2$ and $\cos\alpha$, which allows to insert $\cos \alpha = \eta_2$:
	\begin{align*}
		r(\omega_2, \alpha ) &= \norm{\mathbf{d} - \mathbf{x}} \left( \cos \alpha - \cot \omega_2 \sqrt{1-\cos^2 \alpha} \right)\\
		&= \norm{\mathbf{d} - \mathbf{x}} \left( \eta_2 - \cot \omega_2 \sqrt{1-\left(\eta_2\right)^2} \right) =: r
	\end{align*}
	Inserting these values in the representation of $\mathbf{x}$ yields the intersection point
	\begin{align*}
		\mathbf{y} &= \mathbf{x} + r \begin{pmatrix} \sin ( \omega_1 - \beta + \pi/2)\\ \cos(\omega_1 - \beta + \pi/2)\end{pmatrix}\quad \textrm{if } r > 0
	\end{align*}
	of the circular arcs $\mathfrak{C}(\omega_2, \mathbf{d}, \mathbf{x})$ with the line $\mathfrak{L}(\omega_1, \mathbf{x}, \mathbf{s})$.
\end{proof}

\begin{proof}[Proof of Theorem \ref{thm:secondOrderOperator}]
	Similar to \eqref{eqn:probabilityPhotonScattersOnce}, but inserting a second scattering site and taking into account the physical factors, we can write for the variation of photons scattered twice at points $\mathbf{x}$ and $\mathbf{x}$
	\begin{equation*}
		\mathrm{d}^3N(\mathbf{s},\mathbf{d},\mathbf{x},\mathbf{y}) = I(\mathbf{s}) r_e^4 A_{E_0}(\mathbf{s},\mathbf{x}) A_{E_1}(\mathbf{x}, \mathbf{y}) A_{E_2}(\mathbf{y}, \mathbf{d})  P(\omega_1) P(\omega_2) n_e(\mathbf{x}) n_e(\mathbf{y})\,\mathrm{d}\mathbf{x}\,\mathrm{d}\mathbf{y}\,\mathrm{d}\mathbf{s}.
	\end{equation*}
	This relation is integrated over all possible first scattering sites $\mathbf{x}$. Fixing a measured energy $E$, one further integrates over all possible first scattering angles $\omega_1$. 
	For every $\omega_1$, we can compute the corresponding $\omega_2(\omega_1)$ and the second scattering site $\mathbf{y}(\omega_1) \in \mathfrak{C}(\omega_2, \mathbf{d}, \mathbf{x}) \cap \mathfrak{L}(\omega_1, \mathbf{x}, \mathbf{s})$ by the previous Lemma \ref{lem:SecondScatteringPoint}. 
	We obtain that $g_2$ is proportional to
	\begin{equation*}
		\mathcal{T}_2 (n_e) (E,\mathbf{d},\mathbf{s}) = \int_{\Omega} \int_{(-\pi,\pi)\backslash\{0\}} w_2(\mathbf{x}, \omega_1; \mathbf{d}, \mathbf{s}) n_e(\mathbf{x}) n_e(\mathbf{y}(\omega_1)) \mathrm{d} l(\omega_1)\mathrm{d}\mathbf{x}
	\end{equation*}
	with the weight function $w_2$ as stated in the theorem and the differential line element $\mathrm{d} l(\omega_1)$ encoding the change in the second scattering site $\mathbf{y}$ by the first scattering angle. We remark here that of course $\omega_1 = 0$ is a zero set so it doesn't alter the value of the integral, but it corresponds to the case of no real scattering event and is therefore technically ruled out. 
	
	It remains to compute $\mathrm{d} l(\omega_1) = \norm{\frac{\partial\mathbf{y}}{\partial\omega_1}}\mathrm{d}\omega_1$. For this, we define $\tilde{R}$ to be an orthogonal matrix rotating by an angle $\omega_1 - \beta + \pi/2$ and see that
	\begin{align*}
		\frac{\partial\mathbf{y}}{\partial\omega_1} &= \frac{\partial}{\partial\omega_1} \left( \mathbf{x} + r \begin{pmatrix} \sin ( \omega_1 - \beta + \pi/2)\\ \cos(\omega_1 - \beta + \pi/2)\end{pmatrix} \right)\\
		&= \frac{\partial r}{\partial\omega_1} \begin{pmatrix} \sin ( \omega_1 - \beta + \pi/2)\\ \cos(\omega_1 - \beta + \pi/2)\end{pmatrix} + r \begin{pmatrix} \cos ( \omega_1 - \beta + \pi/2)\\ -\sin(\omega_1 - \beta + \pi/2)\end{pmatrix}\\
		&= \frac{\partial r}{\partial\omega_1} \tilde{R} \begin{pmatrix}0\\ 1\end{pmatrix} + r \tilde{R}\begin{pmatrix}1\\ 0\end{pmatrix} = \tilde{R} \begin{pmatrix}r\\ \frac{\partial r}{\partial\omega_1}\end{pmatrix}
	\end{align*}
	from which we immediately conclude the identity $\mathrm{d} l(\omega_1) = \sqrt{r^2 + \left(\frac{\partial r}{\partial\omega_1}\right)^2}\ \mathrm{d}\omega_1$ where the derivative of $r$ is given by
	\begin{equation*}
		\frac{\partial r}{\partial\omega_1}=\frac{\partial \eta_2}{\partial\omega_1} \left(1+\cot \omega_{2} \frac{\eta_2}{\sqrt{1-\eta_2^2}}\right)-\frac{\sin \omega_{1}}{\sin ^{3} \omega_{2}} \sqrt{1-\eta_2^{2}}
	\end{equation*}
	with $\eta_2$ as in \cref{eqn:etaTwo} and $\partial \eta_2/\partial\omega_1 = R_{1,2} \cos (\omega_1 - \beta + \pi/2) - R_{2,2} \sin(\omega_1 - \beta + \pi/2)$
\end{proof}

\section{\texorpdfstring{Proofs - FIO Properties of $\mathcal{L}_1$}{Proofs - FIO Properties of L1}} \label{pro:ThmL1}
\begin{proof}[Proof of Theorem \ref{thm:fiol1}]
	As $\tilde{n}_e$ is $C^\infty$-smooth, so is the weight $\mathcal{W}_1(\tilde{n}_e)(\mathbf{x}, \mathbf{d}, \mathbf{s})$. Further, the weight does not depend on $\sigma$ at all, so it is immediately clear that it is a symbol of order 0. \\
	Fix the source position $\mathbf{s}$ and parametrize the detector position $\mathbf{d}(\theta)$ by $\theta \in \Theta$ where $\Theta \subset \mathbb{R}$ open. We need to prove that $\Phi(p,\theta,\mathbf{x},\sigma) = - \sigma(p - \phi(\mathbf{x}-\mathbf{s}, \mathbf{d}(\theta)-\mathbf{s}))$ is a non-degenerate phase function.
	Positive homogeneity is clear as $\Phi$ is linear in $\sigma$.
	Further $\partial_{p,\theta}\Phi = -\sigma(\mathrm{d} p - \partial_\theta \phi(\mathbf{x}-\mathbf{s},\mathbf{d}(\theta)-\mathbf{s})) \neq 0$
	for every $(p,\theta,\mathbf{x},\sigma) \in \mathbb{R} \times \Theta\times \Omega \times \mathbb{R}\backslash\{0\}$, in particular $(\partial_{p,\theta}\Phi, \partial_\sigma \Phi)$ doesn't vanish.
	
	Furthermore, it holds $\partial_\mathbf{x} \Phi = \sigma \partial_\mathbf{x} \phi = 0$ if and only if $\nabla_\mathbf{x} \phi = 0$ as $\sigma \neq 0$. For the gradient, we obtain using the definition \eqref{eqn:defphi} of $\phi$
	\begin{multline*}
		\nabla_\mathbf{x} \phi(\mathbf{x}-\mathbf{s},\mathbf{d}-\mathbf{s}) = \frac{1-\rho\kappa}{\norm{\mathbf{x}-\mathbf{s}}(1-\kappa^2)^{3/2}}\left( \frac{\mathbf{d}-\mathbf{s}}{\norm{\mathbf{d}-\mathbf{s}}} - \kappa \frac{\mathbf{x}-\mathbf{s}}{\norm{\mathbf{x}-\mathbf{s}}}\right) - \\
		\frac{\rho}{\norm{\mathbf{x}-\mathbf{s}}(1-\kappa^2)^{1/2}} \frac{\mathbf{x}-\mathbf{s}}{\norm{\mathbf{x}-\mathbf{s}}}.
	\end{multline*}
	The two vectors summed up are orthogonal by the definition of $\kappa$. Hereby, the component in the direction of $(\mathbf{x}-\mathbf{s})/\norm{\mathbf{x}-\mathbf{s}} \neq 0$ is never zero as $\rho \neq 0$, so we have $\nabla_x \phi \neq 0$ and $\phi$ is a phase function.
	
	The phase function is also non-degenerate since we have
	$\partial_\sigma \Phi = (\phi(\mathbf{x}-\mathbf{s},\mathbf{d}(\theta)-\mathbf{s})-p)\mathrm{d}\sigma$
	and therefore $\partial_{\mathbf{x},p,\theta}(\partial_\sigma \Phi) = \left( \partial_\mathbf{x} \phi + \partial_\theta \phi - \mathrm{d} p \right) \mathrm{d} \sigma \neq 0$. Hence $\mathcal{L}_1^{\tilde{n}_e}$ is indeed an FIO.
	It holds $\Omega \subset \mathbb{R}^2$ open and $\Theta \subset \mathbb{R}$ open, i.e. $n = m = 2$. The symbol is of order 0, therefore $\mathcal{L}_1^{\tilde{n}_e}$ is an FIO of order $k = 0 - (2+2-2)/4 = - 1/2$.
\end{proof}

\begin{proof}[Proof of Lemma \ref{lem:immersionl1}]
	Recall the definition of the examined function 
	\begin{equation*}
		\phi(\mathbf{x}-\mathbf{s}, \mathbf{d}-\mathbf{s}) = \frac{\kappa(\mathbf{x}-\mathbf{s},\mathbf{d}-\mathbf{s}) - \rho(\mathbf{x}-\mathbf{s},\mathbf{d}-\mathbf{s})}{\sqrt{1-\kappa(\mathbf{x}-\mathbf{s},\mathbf{d}-\mathbf{s})^2}}.
	\end{equation*}
	Shortening notation, we write $a := \mathbf{x}-\mathbf{s}$, $b := \mathbf{d} - \mathbf{s}$ and further $\kappa := \kappa(a,b)$, $\rho := \rho(a,b)$.
	
	As in the proof of Theorem \ref{thm:fiol1}, we obtain for the gradient
	\begin{equation*}
		\nabla_\mathbf{x} \phi(\mathbf{x}-\mathbf{s},\mathbf{d}-\mathbf{s}) = \frac{1-\rho\kappa}{\norm{a}(1-\kappa^2)^{3/2}}\left( \frac{b}{\norm{b}} - \kappa \frac{a}{\norm{a}}\right) - \frac{\rho}{\norm{a}(1-\kappa^2)^{1/2}} \frac{a}{\norm{a}}.
	\end{equation*}
	At this point, we introduce a change of coordinates by applying the unique orthogonal rotation $R_a$ that maps $a$ to $(\norm{a}, 0)^T$ in $\mathbb{R}^2$. As $R_a$ is a rotation and also independent of $\theta$, when applied to $\nabla_\mathbf{x}\phi$, it does not alter the examined determinant: $\det(\nabla_\mathbf{x}\phi, \partial_\theta \nabla_\mathbf{x}\phi) = \det(R_a\nabla_\mathbf{x}\phi, \partial_\theta R_a \nabla_\mathbf{x}\phi)$. In the new coordinate system, using $R_a b = t(\theta)(\cos(\theta - \xi), \sin(\theta - \xi))^T$ and $\kappa = \cos(\theta - \xi)$ we can simplify the gradient to
	\begin{equation*}
		R_a \nabla_\mathbf{x}\phi(\mathbf{x}-\mathbf{s},\mathbf{d}-\mathbf{s}) = \frac{1}{t(\theta) |\sin(\theta - \xi)|} \begin{pmatrix} -1\\\frac{t(\theta)}{r\sin(\theta-\xi)} - \cot(\theta-\xi) \end{pmatrix}.
	\end{equation*}
	The vector can now be differentiated with respect to theta giving
	\begin{multline*}
		\partial_\theta R_a \nabla_\mathbf{x} \phi = 
		\frac{|\sin(\theta-\xi)|}{t(\theta)\sin^2(\theta-\xi)} \begin{pmatrix}\cot(\theta-\xi) + \frac{t'(\theta)}{t(\theta)}\\ \cot^2(\theta-\xi) + \frac{r-2t(\theta)\cos(\theta-\xi)}{r\sin^2(\theta-\xi)} + \frac{t'(\theta)}{t(\theta)}\cot(\theta-\xi)\end{pmatrix}
	\end{multline*}
	which, after some simplifications, leads to the sought for determinant
	\begin{equation*}
		\det(\nabla_\mathbf{x}\phi, \partial_\theta \nabla_\mathbf{x}\phi) = \frac{1}{rt(\theta)^2\sin^4(\theta-\xi)}(t(\theta)\cos(\theta-\xi) - t'(\theta)\sin(\theta-\xi) - r).
	\end{equation*}
	Assuming condition \eqref{eqn:conditionImmersionl1} is true, the determinant never vanishes, which concludes the proof.
\end{proof}

\section{\texorpdfstring{Proofs - FIO Properties of $\mathcal{L}_2$}{Proofs - FIO Properties of L2}}\label{pro:ThmFioSob}
\begin{proof}[Proof of Theorem \ref{thm:fiol2}]
	The weight function is, as for $\mathcal{L}_1$, a symbol of order zero as it is $C^\infty$-smooth and does not depend on $\sigma$.
	
	The phase function $\Psi$ given by $\Psi(\lambda,\theta,\mathbf{x},\mathbf{y},\sigma) = -\sigma(\lambda - \psi(\mathbf{y}-\mathbf{x},\mathbf{x}-\mathbf{s},\mathbf{d}(\theta)-\mathbf{x}))$ is in $C^\infty((0,2)\times\Theta\times\Omega_2 \times\mathbb{R}\backslash\{0\})$ and obviously homogeneous of order 1 in $\sigma$.
	
	We have to check that $(\partial_{\lambda,\theta}\Psi, \partial_\sigma \Psi)$ and $(\partial_{\mathbf{z}}\Psi, \partial_\sigma \Psi)$ do not vanish for $\sigma \neq 0$. We start with the easier first pair as it holds $\partial_{\lambda,\theta} \Psi = -\sigma(\mathrm{d}\lambda - \partial_\theta \psi(\mathbf{x}-\mathbf{x},\mathbf{x}-\mathbf{s},\mathbf{d}(\theta)-\mathbf{x})) \neq 0$
	when $\sigma \neq 0$. Thus $(\partial_{\lambda,\theta}\Psi, \partial_\sigma \Psi)$ does not vanish. 
	
	The first two components of $\partial_\mathbf{z} \Psi$ are $\partial_\mathbf{x} \Psi = \sigma \partial_\mathbf{x} \psi(\mathbf{y}-\mathbf{x},\mathbf{x}-\mathbf{s},\mathbf{d}(\theta)-\mathbf{x})$
	which is zero if and only if $\nabla_\mathbf{x}\psi(\mathbf{y}-\mathbf{x},\mathbf{x}-\mathbf{s},\mathbf{d}(\theta)-\mathbf{x}) = 0$. 
	Abbreviating $a := \mathbf{y}-\mathbf{x}, b := \mathbf{x}-\mathbf{s}$ and $c := \mathbf{d} - \mathbf{x}$, where $a,b,c \neq 0$, the gradient is given by
	\begin{equation}
		\nabla_\mathbf{x} (\psi(\mathbf{y}-\mathbf{x},\mathbf{x}-\mathbf{s},\mathbf{d}-\mathbf{x})) = (-\nabla_a + \nabla_b) \kappa(a,b) - \frac{1}{(1+\phi(a,c)^2)^{3/2}}(\nabla_a + \nabla_c)\phi(a,c)\label{eqn:gradxPsi}.
	\end{equation}
	By straightforward differentiation by $a$ and the symmetry $\kappa(a,b) =\kappa(b,a)$ we obtain
	\begin{equation*}
		\nabla_a\kappa(a,b) = \frac{b}{\norm{a}\norm{b}} - \frac{a\cdot b}{\norm{a}^3 \norm{b}} a, \quad \nabla_b\kappa(a,b) = \frac{a}{\norm{a}\norm{b}} - \frac{a\cdot b}{\norm{a} \norm{b}^3} b
	\end{equation*}
	and, writing $\bar{\rho} = \rho(a,b) = \norm{a} / \norm{b}$ as well as $\bar{\kappa} = \kappa(a,b)$, their difference can be simplified to
	\begin{equation}\label{eqn:gradxPsiI}
		(-\nabla_a + \nabla_b) \kappa(a,b) = \left(\bar{\kappa} + \bar{\rho}\right) \frac{a}{\norm{a}^2} - \left(\bar{\kappa} + \bar{\rho}^{-1}\right) \frac{b}{\norm{b}^2}.
	\end{equation}
	Next, we compute the terms depending on $(a,c)$. Abbreviating $\kappa = \kappa(a,c)$, $\rho = \rho(a,c)$, we get
	\begin{align*}
		(1+\phi(a,c)^2)^{-3/2} &= \left(1 + \left(\frac{\kappa-\rho}{1-\kappa^2}\right)^2\right)^{3/2} = \left( \frac{1-\kappa^2}{1-2\kappa\rho + \rho^2}\right)^{3/2},\\
		\nabla_a \phi(a,c) &= \frac{1}{\norm{a}(1-\kappa^2)^{3/2}} \left( (-\kappa -\rho + 2\rho\kappa^2) \frac{a}{\norm{a}} + (1-\rho\kappa) \frac{c}{\norm{c}} \right),\\
		\nabla_c \phi(a,c) &= \frac{1}{\norm{c}(1-\kappa^2)^{3/2}} \left( (1-\rho\kappa) \frac{a}{\norm{a}} - (\kappa - \rho) \frac{c}{\norm{c}} \right).
	\end{align*}
	The combination of these three terms simplifies then to
	\begin{equation}\label{eqn:gradxPsiII}
		- \frac{1}{(1+\phi(a,c)^2)^{3/2}}(\nabla_a + \nabla_c)\phi(a,c) = - \frac{1}{\norm{a}(1-2\kappa\rho+\rho^2)^{1/2}} \left( -\kappa\frac{a}{\norm{a}} + \frac{c}{\norm{c}} \right)
	\end{equation}
	It follows that, by $\kappa = \frac{a\cdot c}{\norm{a}\norm{c}}$, \eqref{eqn:gradxPsiII} is orthogonal to $a$. This shows that it is sufficient to prove that eq. \eqref{eqn:gradxPsiI} always has a nonzero component in the direction of $a$ in order to prove that the sum of the two terms never vanishes. We therefore compute 
	\begin{equation*}
		((-\nabla_a + \nabla_b) \kappa(a,b))^T a = \left(\bar{\kappa} + \bar{\rho}\right) \frac{a^T a}{\norm{a}^2} - \left(\bar{\kappa} + \bar{\rho}^{-1}\right) \frac{b^T a}{\norm{b}^2} = \frac{\norm{a}^2\norm{b}^2 - (a^T b)^2}{\norm{a}\norm{b}^3}.
	\end{equation*}
	The last term obviously vanishes if and only if $a = \mathbf{y}-\mathbf{x}$ and $b = \mathbf{x}-\mathbf{s}$ are collinear, a case that was discarded in the definition of $\Omega_2$. Therefore, the component of $(-\nabla_a + \nabla_b) \kappa(a,b)$ in the direction of $a$ is nonzero for $\mathbf{z} \in \Omega_2$. We have proved that the gradient $\nabla_\mathbf{x} \psi$ and thus also $\partial_\mathbf{z} \Psi$ do not vanish, therefore $\Psi$ is a phase function.
	
	It remains to prove that the phase in non-degenerate. The necessary condition is always satisfied as $\partial_{\mathbf{z},\lambda,\theta} (\partial_\sigma \Psi) = \left( \partial_\mathbf{z} \psi + \partial_\theta \psi - \mathrm{d} \lambda \right) \mathrm{d} \sigma \neq 0$.
	As $\Omega_2 \subset \mathbb{R}^4$ open and $\Theta \subset \mathbb{R}$ open, it holds $n=4$, $m = 2$ and $\mathcal{L}_2^{\tilde{n}_e}$ is an FIO of order $\frac{1}{2} - \frac{4+2}{4} = -1$.
\end{proof}

\begin{proof}[Proof of Theorem \ref{thm:sobolevtheoreml2}]
	In view of Theorem \ref{thm:fiol2} and Lemma \ref{lem:sobolevMappingFIO}, it remains to show that the vectors $(\nabla_\mathbf{x} \psi, \nabla_\mathbf{y} \psi)^T, (\partial_\theta \nabla_\mathbf{x} \psi, \partial_\theta \nabla_\mathbf{y} \psi)^T$
	are linearly independent for every $(\mathbf{d}(\theta),\mathbf{x},\mathbf{y})$, $\theta \in \Theta$, $(\mathbf{x},\mathbf{y}) \in \Omega_2$.
	
	We can recycle some of the computations in the proof of Theorem \ref{thm:fiol2}. It suffices to show that already the first two components of the gradients are linearly independent to prove the Theorem. 
	
	The first two components of the first vector, $\nabla_\mathbf{x} \psi$, were computed in eq. \eqref{eqn:gradxPsi} and split up into two parts: $\nabla_\mathbf{x} \psi = r_1(a,b) + r_2(a,c)$,
	where $r_1(a,b)$ given in eq. \eqref{eqn:gradxPsiI}, depending only on $a = \mathbf{y}-\mathbf{x}$ and $b=\mathbf{x}-\mathbf{s}$, and $r_2(a,c)$ given in eq. \eqref{eqn:gradxPsiII}, depending only on $a$ and $c=\mathbf{d}(\theta)-\mathbf{x}$. We then argued that $r_2(a,c)^T a = 0$.
	
	Consider now the second vector $\partial_\theta \nabla_\mathbf{x} \psi = \partial_\theta r_1(a,b) + \partial_\theta r_2(a,c)$.
	$r_1(a,b)$ does not depend on $\theta$ at all, so $\partial_\theta r_1(a,b) = 0$. The second part $r_2(a,c)$ is orthogonal to $a$ and so is its derivative $\partial_\theta r_2(a,c)$.
	The vector $\partial_\theta \nabla_\mathbf{x} \psi$ is thus a multiple of $a^\perp$, and so is $r_2(a,c)$. In order to prove that $\partial_\theta \nabla_\mathbf{x} \psi$ and $\nabla_\mathbf{x} \psi$ are linearly independent, it is therefore sufficient to show that $r_1(a,b)$ has a nonzero component in the direction of $a$. But this leads us back to the proof of Theorem \ref{thm:fiol2}, where we studied this exact case and showed that it is always satisfied when $\mathbf{y}-\mathbf{x}$ and $\mathbf{x}-\mathbf{s}$ are not collinear, namely when $(\mathbf{x},\mathbf{y})\in \Omega_2$. This proves the Theorem.
\end{proof}

\section*{Acknowledgments} The second author was supported by the Deutsche Forschungsgemeinschaft (DFG) under the grant RI 2772/2-1 and the Cluster of Excellence EXC 2075 "Data-Integrated Simulation Science" at University of Stuttgart. 

\printbibliography

@article{Rigaud21,
author={G. Rigaud},
title={{3D Compton scattering imaging with multiple scattering: analysis by FIO and contour reconstruction}},
journal={Inverse Problems},
volume = {37},
number = {6},
doi = {10.1088/1361-6420/abf22b},
year={2021},
}

@ARTICLE{Wang2006,
	author={Jing Wang and Tianfang Li and Hongbing Lu and Zhengrong Liang},
	journal={IEEE Transactions on Medical Imaging}, 
	title={Penalized weighted least-squares approach to sinogram noise reduction and image reconstruction for low-dose X-ray computed tomography}, 
	year={2006},
	volume={25},
	number={10},
	pages={1272-1283},
	doi={10.1109/TMI.2006.882141}}

@article{PhysRevA.13.335,
  title = {{Multiple scattering in the Compton effect. I. Analytic treatment of angular distributions and total scattering probabilities}},
  author = {A.C. Tanner and I.R. Epstein},
  journal = {Phys. Rev. A},
  volume = {13},
  issue = {1},
  pages = {335--348},
  numpages = {0},
  year = {1976},
  month = 1,
  publisher = {American Physical Society},
  doi = {10.1103/PhysRevA.13.335},
}

@article{PhysRevA.14.313,
  title = {{Multiple scattering in the Compton effect. II. Analytic and numerical treatment of energy profiles}},
  author = {A.C. Tanner and I.R. Epstein},
  journal = {Phys. Rev. A},
  volume = {14},
  issue = {1},
  pages = {313--327},
  numpages = {0},
  year = {1976},
  month = 7,
  publisher = {American Physical Society},
  doi = {10.1103/PhysRevA.14.313},
}

@article{PhysRevA.14.328,
  title = {{Multiple scattering in the Compton effect. III. Monte Carlo calculations}},
  author =  {A.C. Tanner and I.R. Epstein},
  journal = {Phys. Rev. A},
  volume = {14},
  issue = {1},
  pages = {328--340},
  numpages = {0},
  year = {1976},
  month = 7,
  publisher = {American Physical Society},
  doi = {10.1103/PhysRevA.14.328}
}

@article{Stonestrom81,  
author={J. P. {Stonestrom} and R. E. {Alvarez} and A. {Macovski}},
journal={IEEE Transactions on Biomedical Engineering},   
title={A Framework for Spectral Artifact Corrections in X-Ray CT},
year={1981},  
volume={BME-28},  
number={2},  
pages={128-141},
doi = {10.1109/TBME.1981.324786}
}

@book{Leroy11,
author = {Claude {Leroy} and Pier-Giorgio {Rancoita}},
title = {Principles of radiation interaction in matter and detection},
doi = {10.1142/9167},
publisher = {World Scientific},
address = {Singapore},
year = {2011},
}

@article{Sidky08,
year = 2008,
publisher = {{IOP} Publishing},
volume = {53},
number = {17},
pages = {4777--4807},
author = {Emil Y Sidky and Xiaochuan Pan},
title = {Image reconstruction in circular cone-beam computed tomography by constrained, total-variation minimization},
doi = {10.1088/0031-9155/53/17/021},
journal = {Physics in Medicine and Biology},
}

@article{Almansa08,
Author = {Almansa, A. and Ballester, C. and Caselles, V. and Haro, G.},
Journal = {Journal of Scientific Computing},
Number = {3},
Pages = {209--236},
doi={10.1007/s10915-007-9160-x},
Title = {A TV Based Restoration Model with Local Constraints},
Volume = {34},
Year = {2008},
}

@article{Zhu13,
Author = {Zhu, Zangen and Wahid, Khan and Babyn, Paul and Cooper, David and Pratt, Isaac and Carter, Yasmin},
Journal = {Computational and Mathematical Methods in Medicine},
Pages = {185750},
Title = {Improved Compressed Sensing-Based Algorithm for Sparse-View CT Image Reconstruction},
Volume = {2013},
doi={10.1155/2013/185750},
Year = {2013}
}

@article{Klein29,
Author = {Oskar Klein and Yoshio Nishina},
Journal = {Zeitschrift für Physik},
Number = {11},
Pages = {853--868},
Title = "{Über die Streuung von Strahlung durch freie Elektronen nach der neuen relativistischen Quantendynamik von Dirac}",
Volume = {52},
doi={10.1007/BF01366453},
Year = {1929}
}

@Inbook{Krishnan15,
author="Krishnan, Venkateswaran P.
and Quinto, Eric Todd",
editor="Scherzer, Otmar",
doi={},
title="Microlocal Analysis in Tomography",
bookTitle="Handbook of Mathematical Methods in Imaging",
year="2015",
publisher="Springer New York",
address="New York, NY",
pages="847--902",
}

@article{WebberIPI19,
title = {{Microlocal analysis of a spindle transform}},
journal = {{Inverse Problems \& Imaging}},
volume = {13},
number = {2},
pages = {231-261},
year = {2019},
author = {J. W. Webber and S. Holman},
}

@misc{Webber19,
author = {{Webber}, James and {Quinto}, Eric Todd},
title = "{Microlocal analysis of a Compton tomography problem}",
year = 2019,
archivePrefix = {arXiv},
eprint = {1902.09623},
eprintclass = {math.FA},
}

@article{Rudin92,
author = {Rudin, Leonid I. and Osher, Stanley and Fatemi, Emad},
title = {Nonlinear Total Variation Based Noise Removal Algorithms},
year = {1992},
publisher = {Elsevier Science Publishers B. V.},
address = {NLD},
volume = {60},
number = {1–4},
journal = {Phys. D},
pages = {259–268},
doi = {10.1016/0167-2789(92)90242-F},
}

@article{Acar94,
year = 1994,
publisher = {{IOP} Publishing},
volume = {10},
number = {6},
pages = {1217--1229},
author = {R Acar and C R Vogel},
title = {Analysis of bounded variation penalty methods for ill-posed problems},
doi = {10.1088/0266-5611/10/6/003},
journal = {Inverse Problems},
}

@article{Rigaud18,
year = 2018,
publisher = {{IOP} Publishing},
volume = {34},
number = {7},
pages = {075004},
author = {Gaël Rigaud and Bernadette N Hahn},
title = {3D Compton scattering imaging and contour reconstruction for a class of Radon transforms},
doi={10.1088/1361-6420/aabf0b},
journal = {Inverse Problems}
}

@article{Cebeiro21,
	year = 2021,
	month = 4,
	publisher = {{IOP} Publishing},
	volume = {37},
	number = {5},
	pages = {054001},
	author = {J Cebeiro and C Tarpau and M A Morvidone and D Rubio and M K Nguyen},
	title = {On a three-dimensional Compton scattering tomography system with fixed source},
	journal = {Inverse Problems},
}

@article{Hoermander71,
author = "Hörmander, Lars",
fjournal = "Acta Mathematica",
journal = "Acta Math.",
pages = "79--183",
publisher = "Institut Mittag-Leffler",
title = "Fourier integral operators. I",
volume = "127",
doi={10.1007/BF02392052},
year = "1971"
}

@article{Kanematsu16,
year = 2016,
publisher = {{IOP} Publishing},
volume = {61},
number = {13},
pages = {5037--5050},
doi={10.1088/0031-9155/61/13/5037},
author = {Nobuyuki Kanematsu and Taku Inaniwa and Minoru Nakao},
title = {Modeling of body tissues for Monte Carlo simulation of radiotherapy treatments planned with conventional x-ray {CT} systems},
journal = {Physics in Medicine and Biology},
}

@article{Webber18,
year = 2018,
publisher = {{IOP} Publishing},
volume = {34},
number = {8},
pages = {084001},
author = {James W Webber and William R B Lionheart},
title = {Three dimensional Compton scattering tomography},
doi = {10.1088/1361-6420/aac51e},
journal = {Inverse Problems},
}

@article{Webber20,
year = 2020,
publisher = {{IOP} Publishing},
volume = {36},
number = {2},
pages = {025007},
author = {James Webber and Eric L Miller},
title = {Compton scattering tomography in translational geometries},
doi = {10.1088/1361-6420/ab4a32},
journal = {Inverse Problems},
}

@misc{Goedecke2022,
	title={Imaging based on Compton scattering: model-uncertainty and data-driven reconstruction methods}, 
	author={Janek Gödeke and Gaël Rigaud},
	year={2022},
	eprint={2202.00810},
	archivePrefix={arXiv},
	primaryClass={math.NA}
}

@article{Wang99,
author = {Wang, J.  and Chi, Z. and Wang, Y.},
title = {Analytic reconstruction of Compton scattering tomography},
doi={10.1063/1.370949},
journal = {Journal of Applied Physics},
volume = {86},
number = {3},
pages = {1693-1698},
year = {1999},
}

@inbook{Truong12,
year = 2012,
editor = {Andriychuk, M.},
title = {Recent Developments on Compton Scatter Tomography: Theory and Numerical Simulations},
doi={10.5772/50012},
author = {Truong, T. T. and Nguyen, M. K.},
bookTitle = {Numerical Simulation - From Theory to Industry},
publisher = {Intech},
}

@article{Truong07,
Author = {Truong, T. T. and Nguyen, M. K. and Zaidi, H.},
Journal = {International Journal of Biomedical Imaging},
Pages = {092780},
doi={10.1155/2007/92780},
Title = {The Mathematical Foundations of 3D Compton Scatter Emission Imaging},
Volume = {2007},
Year = {2007}
}

@misc{Tarpau20a,
author = {Tarpau, C. and Cebeiro, J. and Nguyen, M. K. and
Rollet, G. and Dumas, L.},
title = {On the design of a CST system and its extension to a bi-imaging modality},
year = 2020,
archivePrefix = {arXiv},
eprint = {2007.02750},
eprintclass = {physics.med-ph},
}

@article{Shrimpton81,
year = 1981,
publisher = {{IOP} Publishing},
volume = {26},
number = {5},
pages = {907--911},
author = {Shrimpton, P. C.},
title = {Electron density values of various human tissues: in vitro Compton scatter measurements and calculated ranges},
	doi = {10.1088/0031-9155/26/5/010},
journal = {Physics in Medicine and Biology},
}

@article{Lale59,
year = 1959,
publisher = {{IOP} Publishing},
volume = {4},
number = {2},
pages = {159--167},
author = {Lale, P. G.},
title = {The Examination of Internal Tissues, using Gamma-ray Scatter with a Possible Extension to Megavoltage Radiography},
journal = {Physics in Medicine and Biology},
doi={10.1088/0031-9155/4/2/305},
}

@book{Giusti84,
author = {Giusti, E.},
title = {Minimal surfaces and functions of bounded variation},
publisher = {Springer Boston},
doi={10.1007/978-1-4684-9486-0},
year = {1984},
}

@article{Clarke73,
year = 1973,
publisher = {{IOP} Publishing},
volume = {18},
number = {4},
pages = {532--539},
author = {Clarke, R. L. and Van Dyk, G.},
title = {A new method for measurement of bone mineral content using both transmitted and scattered beams of gamma-rays},
doi={10.1088/0031-9155/18/4/005},
journal = {Physics in Medicine and Biology},
}

@article{Brateman84,
year = 1984,
publisher = {{IOP} Publishing},
volume = {29},
number = {11},
pages = {1353--1370},
doi = {10.1088/0031-9155/29/11/004},
author = {Brateman, L. and Jacobs, A. M. and Fitzgerald, L. T.},
title = {Compton scatter axial tomography with x-rays: {SCAT}-{CAT}},
journal = {Physics in Medicine and Biology},
}

@misc{xcomDatabase,
author = {Berger, M. J. and Hubbell, J. H. and Seltzer, S. M. and Chang, J. and Coursey, J. S. and Sukumar, R. and Zucker, D. S. and Olsen, K.},
title = {NIST XCOM: Photon Cross Sections Database},
note = {Accessed: 2020-09-05}, 
url = {http://physics.nist.gov/xcom},
}

@incollection{Burger13,
author = {Burger, Martin and Osher, Stanley},
booktitle = {Level Set and PDE-based Reconstruction Methods},
pages = {1-70},
peerreviewed = {Yes},
series = {Lecture Notes in Mathematics},
title = {{A} guide to the {TV} zoo},
year = {2013}
}

@article{Webber15,
title={X-ray Compton scattering tomography},
author={Webber, James},
journal={Inverse Problems in Science and Engineering},
year={2015},
volume={24},
pages={1323 - 1346},
doi={10.1080/17415977.2015.1104307}
}

@article{Norton94,
author = {Norton, Stephen J.},
title = {Compton scattering tomography},
journal = {Journal of Applied Physics},
doi={10.1063/1.357668},
volume = {76},
number = {4},
pages = {2007-2015},
year = {1994},
}

@article{Nguyen10,
year = 2010,
publisher = {{IOP} Publishing},
volume = {26},
number = {6},
pages = {065005},
author = {M K Nguyen and T T Truong},
doi = {10.1088/0266-5611/26/6/065005},
title = {Inversion of a new circular-arc Radon transform for Compton scattering tomography},
journal = {Inverse Problems},
}

@article{Palamodov11,
year = 2011,
publisher = {{IOP} Publishing},
volume = {27},
number = {12},
pages = {125004},
author = {V P Palamodov},
title = {An analytic reconstruction for the Compton scattering tomography in a plane},
doi = {10.1088/0266-5611/27/12/125004},
journal = {Inverse Problems},
}

@book{Natterer01A,
author = {Natterer, F.},
title = {The Mathematics of Computerized Tomography},
doi={10.1137/1.9780898719284},
publisher = {Society for Industrial and Applied Mathematics},
year = {2001},
}

@book{Natterer01B,
author = {Natterer, Frank and Wübbeling, Frank},
title = {Mathematical Methods in Image Reconstruction},
doi={10.1137/1.9780898718324},
publisher = {Society for Industrial and Applied Mathematics},
year = {2001},
}

@INPROCEEDINGS{Evans97,
author={B. L. {Evans} and J. B. {Martin} and L. W. {Burggraf} and M. C. {Roggemann}},
booktitle={1997 IEEE Nuclear Science Symposium Conference Record}, 
title={Nondestructive inspection using Compton scatter tomography}, 
doi = {10.1109/NSSMIC.1997.672608},
year={1997},
volume={1},
pages={386-390},
}

@article{Rigaud17,
author = {Rigaud, Gaël},
year = {2017},
pages = {2217-2249},
title = {Compton Scattering Tomography: Feature Reconstruction and Rotation-Free Modality},
doi={10.1137/17M1120105},
volume = {10},
journal = {SIAM Journal on Imaging Sciences},
}

@article{Hahn19,
author = {Hahn, B. N. and Kienle Garrido, M.-L.},
journal = {Inverse Problems},
number = 9,
pages = 094005,
title = {An efficient reconstruction approach for a class of dynamic imaging operators},
doi = {10.1088/1361-6420/ab178b},
volume = 35,
year = 2019
}

@article{Hahn14,
author = {Hahn, B. N.},
journal = {Journal of Inverse and Ill-posed Problems},
number = 3,
pages = {323-339},
title = {Reconstruction of dynamic objects with affine deformations in computerized tomography},
doi={10.1515/jip-2012-0094},
volume = 22,
year = 2014
}

@article{Kuchment95,
year = 1995,
publisher = {{IOP} Publishing},
volume = {11},
number = {3},
pages = {571--589},
author = {P Kuchment and K Lancaster and L Mogilevskaya},
title = {On local tomography},
doi = {10.1088/0266-5611/11/3/006},
journal = {Inverse Problems},
}

@article{Greenleaf2002,
  title={Oscillatory and Fourier integral operators with degenerate canonical relations},
  author={A. Greenleaf and A. Seeger},
  journal={Publicacions Matematiques},
  year={2002},
  volume={46},
  pages={93-141}
}

\end{document}